\documentclass[12pt
]{amsart}

\usepackage[margin=1in]{geometry}

\usepackage{amsmath}
\usepackage{amsfonts}
\usepackage{amssymb}
\usepackage{amsthm}
\usepackage{setspace}
\usepackage{amsrefs}
\usepackage{hyperref}
\hypersetup{hidelinks}

\title[H\"{o}lder estimates for homotopy operators]{H\"{o}lder estimates for homotopy operators on strictly pseudoconvex domains with $C^2$ boundary}
\author[]{Xianghong Gong}

 \address{Department of Mathematics,
 University of Wisconsin-Madison, Madison, WI 53706, U.S.A.}
 \email{gong@math.wisc.edu}


\newcommand{\dist}{\operatorname{dist}}

\newtheorem{thm}{Theorem}[section]
\newtheorem{cor}[thm]{Corollary}
\newtheorem{prop}[thm]{Proposition}
\newtheorem{lemma}[thm]{Lemma}

\theoremstyle{definition}

\newtheorem{defn}[thm]{Definition}
\newtheorem{exmp}[thm]{Example}

\newtheorem{rem}[thm]{Remark}

\renewcommand{\th}[1]{\begin{thm}\label{#1}}
\renewcommand{\eth}{\end{thm}}
\newcommand{\co}[1]{\begin{cor}\label{#1}}
\newcommand{\eco}{\end{cor}}
\renewcommand{\le}[1]{\begin{lemma}\label{#1}}
\newcommand{\ele}{\end{lemma}}
\newcommand{\pr}[1]{\begin{prop}\label{#1}}
\newcommand{\epr}{\end{prop}}

\newcommand{\ga}{\begin{gather}}
\newcommand{\ega}{\end{gather}}
\newcommand{\gan}{\begin{gather*}}
\newcommand{\egan}{\end{gather*}}
\newcommand{\al}{\begin{align}}
\newcommand{\eal}{\end{align}}
\newcommand{\aln}{\begin{align*}}
\newcommand{\ealn}{\end{align*}}
\newcommand{\eq}[1]{\begin{equation}\label{#1}}
\newcommand{\eeq}{\end{equation}}

\newcommand{\DD}[2]{\frac{\partial #1}{\partial #2}}
\newcommand{\f}[2]{\frac{#1}{#2}}

\newcommand{\ci}{~\cite}

\newcommand{\cc}{{\bf C}}
\newcommand{\nn}{{\bf N}}

\newcommand{\rr}{{\bf R}}

\newcommand{\oo}{\"{o}}

\newcommand{\ov}{\overline}

\newcommand{\RE}{\operatorname{Re}}
\newcommand{\IM}{\operatorname{Im}}
\renewcommand{\dbar}{\overline\partial}

\newcommand{\diam}{\operatorname{diam}}

\newcommand{\cL}{\mathcal}

\newcommand{\all}{\alpha}

\newcommand{\gaa}{\gamma}

\newcommand{\del}{\delta}
\newcommand{\Del}{\Delta}
\newcommand{\var}{\varphi}
\newcommand{\e}{\epsilon}
\newcommand{\om}{\omega}
\newcommand{\Om}{\Omega}

\newcommand{\la}{\lambda}

\newcommand{\pd}{\partial}

\newcommand{\yt}{\frac{1}{2}}

\newcommand{\re}[1]{(\ref{#1})}
\newcommand{\rea}[1]{$(\ref{#1})$}
\newcommand{\rl}[1]{Lemma~\ref{#1}}
\newcommand{\nrc}[1]{Corollary~\ref{#1}}
\newcommand{\rp}[1]{Proposition~\ref{#1}}
\newcommand{\rt}[1]{Theorem~\ref{#1}}
\newcommand{\rd}[1]{Definition~\ref{#1}}

\newcommand{\rla}[1]{Lemma~$\ref{#1}$}

\newcommand{\rpa}[1]{Proposition~$\ref{#1}$}

\newcommand{\rda}[1]{Definition~$\ref{#1}$}

\newcommand{\supp}{\operatorname{supp}}

\newcommand{\db}{\dbar}
\newcommand{\dbb}{\dbar_b}

\newcounter{pp}
\newcommand{\bpp}{\begin{list}{$\hspace{-1em}\alph{pp})$}{\usecounter{pp}}}
\newcommand{\epp}{\end{list}}

\newcounter{ppp}
\newcommand{\bppp}{\begin{list}{$\hspace{-1em}(\roman{ppp})$}{\usecounter{ppp}}}
\newcommand{\eppp}{\end{list}}

\def\beq{\begin{equation}}
\def\eeq{\end{equation}}

\begin{document}
\begin{abstract}
We  derive a new homotopy formula for a strictly pseudoconvex domain of $C^2$ boundary in ${\mathbf C}^n$ by using a method of Lieb and Range and obtain estimates in Lipschitz spaces for the homotopy operators. For $r>1$ and $q>0$, we obtain a $\Lambda_{r+{1}/{2}}$ solution $u$  to $\overline\partial u=f$ for a $\overline\partial$-closed $(0,q)$-form $f$ of class $\Lambda_{r}$  in the domain.  We apply the estimates   to obtain  boundary regularities of $\mathcal D$-solutions for a domain in $\cc^n\times\rr^m$.  
\end{abstract}

 \maketitle


\setcounter{thm}{0}\setcounter{equation}{0}
\section{Introduction}\label{sec1}

The main purpose of this paper is to  show the boundary regularity for $\db$ solutions in a strictly pseudoconvex domain $D$ in $\cc^n$ under the minimal smoothness condition of the boundary $\pd D\in C^2$. We will also derive a homotopy formula for the domain $D$,
\eq{va0qH}
\var=\db H_q\var+H_{q+1}\db\var, \quad q\geq1
\eeq
that admits a derivative estimate. Here $\var $ is a $(0,q)$-form in $\ov D$  and $\var$, $\db\var$
   are in $C^1(\ov D)$. We will prove the following $C^{r+1/2}$ estimate.
\th{thm1} Let $r\in[1,\infty)$ and $1\leq q\leq n$.
Let $D$ be a bounded strictly pseudoconvex domain  of  $C^{2}$ boundary in $\cc^n$.  If $r+1/2$ is  non integer, then
\eq{hqpsi}
|H_q\var|_{C^{r+1/2}(\ov D)}\leq C_r(D)|\var|_{C^r(\ov D)},
\eeq
where $C_r(D)<\infty$ depends only on $r$ and the domain $D$.
\eth

The study of regularities of $\dbar$ solutions via integral representations  has a long history.
The sup-norm estimate of $\db$ solutions  was proved by  Grauert-Lieb~\cite{GL70} and  Henkin\ci{He69} for
$(0,1)$-forms (the forms are thus $\db$-closed).   Kerzman~\ci{Ke71} obtained  $L^p$ and $C^{\beta}$ estimates of $\db$ solutions for $(0,1)$-forms and all $\beta<1/2$, and {\O}vrelid~\ci{Ov71} obtained a homotopy formula with  homotopy operators admitting $L^p$ estimates for all $(0,q)$-forms.
  Lieb~\ci{Li70} obtained the $L^\infty$ and the $C^\beta$ estimates of $\db$ solutions for $(0,q)$-forms. Finally,
Henkin and Romanov~\ci{RH71}   achieved the $C^{1/2}$ estimate of $\db$ solutions for continuous $(0,1)$-forms. The two $\db$-solution operators in~\cites{He69,GL70} make essential uses of the Henkin-Ram{\'{\i}}rez functions constructed indepently by Henkin \ci{He69} and Ram{\'{\i}}rez\ci{Ra69}.
On the other hand,
Stein showed that the $C^{1/2}$ estimate
 is optimal for  $\db$-closed continuous $(0,1)$-forms in the unit ball of $\cc^n$ for $n>1$ (see~\cites{Ke71,HL84}).
Note that Treves~\ci{Tr92}   studied the boundary regularity for the Leray-Koppelman homotopy operator.
Noticeably   
 the $C^{1/2}$  estimate, valid for for all continuous   $(0,q)$-forms that are not necessarily $\db$-closed,  was first obtained    by Range-Siu~\ci{RS73} for a homotopy operator $T_q$.
 However, to the author's best knowledge it remains open if   the 
 $\db$ solution operators in the above-mentioned results   have a boundary regularity beyond the $C^{1/2}$ estimate when they act on the continuous forms that are not $\db$-closed.
  There are, of course, important results under the conditions that   $\var$ is $\db$-closed
  and $r$ is a positive integer $k$:
 Siu~\ci{Si74}  proved the $C^{k+1/2}$ estimate   for $T_1$ and Alt~\ci{Al74} obtained analogous results for the two $\db$ solution operators of  Kerzman~\ci{Ke71} and  Grauert-Lieb~\ci{GL70} for $(0,1)$-forms.  For  $\db$-closed $(0,q)$-forms $\var$ with $q\geq1$,  Lieb and Range~\ci{LR80} constructed a  $\db$ solution operator $H_q$  and proved  \re{hqpsi}  when $\pd D\in C^{k+2}$, and in~\cites{LR86, LR86a} they also showed that  Kohn's canonical solution  $u$ to $\db u=\var$ is in $C^{k+1/2}(\ov D)$ when  $  \var\in C^{k}(\ov D)$  and   $\pd D\in C^\infty$. The above-mentioned results are for   strictly pseudoconvex domains.  Range~\ci{Ra90} obtained a H\"older estimate for $\db$ solutions in finite type pseudoconvex domains of $\cc^2$. There are derivative estimates for $\dbar$ solutions in    convex domains of D'Angelo finite type $m$:
  Diederich-Forn\ae ss-Wiegerlinck~\ci{DFW86} obtained the $C^{1/m}$ estimate for ellipsoids,  Diederich-Fischer-Forn\ae ss~\ci{DFF99}  and Cumenge~\ci{Cu01} obtained  the $C^{1/m}$ estimate,  and Alexandre~\ci{Al06} achieved the $C^{k+1/m}$ estimate for $\db$ solutions.

  \rt{thm1} does not require that $\var$ is $\db$-closed. Here are some related results.
 An interior estimate of gaining one derivative  for $\var\in C^r$ with a non-integer $r$ was obtained by Webster~\ci{We89}.  We  mention that  Alexandre~\ci{Al06} obtained  the $C^{1/m}$ estimate for a homotopy operator   on convex finite-type domains.  For the $\db_b$ operator in a domain
  in a strictly pseudoconvex hypersurface $M$ in $\cc^n$ with $n\geq4$, the interior $C^k$ estimate was obtained by  Webster~\ci{We89cr} and Ma-Michel~\cite{MM93} proved a boundary regularity for   homotopy operators. Gong-Webster~\ci{GW11} obtained an interior $C^{k+1/2}$ estimate for Henkin homotopy operators
  when the $M$ is in $ C^{k+2}$.
   Range and Siu~\ci{RS73} proved the $C^\beta$ estimate for all $\beta<1/2$ for $\db$ solutions of continuous $(0,q)$-forms on the transversal intersection of
  strictly pseudoconvex domains; see also Poljakov~\ci{Po72} for   related results.
 It is open if the $C^{1/2}$   estimate  holds $\db$ solutions for continuous forms in this situation.
   Higher order derivative estimates for $\db$ solutions were obtained by Brinkmann~\ci{Br84}, Michel~\ci{Ma88}, and Michel-Perotti~\ci{MP90} for the intersection.
 Peters~\ci{Pe91}  constructed a new homotopy operator for the weakly transversal intersection of strictly pseudoconvex domains and obtained higher order derivative estimates with some loss of derivatives.  Note that all these results require the boundary of domains to be sufficiently smooth. For weakly pseudoconvex domains with $C^\infty$ boundary,  Michel~\ci{Mi91}  and Michel-Shaw~\ci{MS99} respectively constructed  homotopy operators with the $C^\infty$ regularity for the domains and for their transversal intersection.

We will derive a   homotopy operator for a strictly pseudoconvex   domain $D$ with $C^2$ boundary, by perfecting the formulation of the  Lieb-Range $\db$ solution operator.   The   homotopy operator 
has the form
\ga\label{nhq}
H_q\var(z)=\int_{\cc^n}\Om_{0,q-1}^0(z,\zeta)\wedge E\var(\zeta)+\int_{\cc^n\setminus D}\Om_{0,q-1}^{01}(z,\zeta)\wedge[\db,E]\var(\zeta)
\end{gather}
for $z\in D$ and $q>0$. Here    $E\colon C(\ov D)\to C_0(\cc^n)$ is a linear extension operator constructed by Calder\'on~\ci{Ca61} and Stein~\ci{St70}, and it satisfies two important properties
$$ 
|Ef |_{\cc^n;r}\leq C_r |f|_{\ov D;r},\quad |Ef |_{\Lambda_r(\cc^n)}\leq C_r |f|_{\Lambda_r(\ov D)}.
$$ 
 Here $C^r(\ov D)$  with norm $|\cdot|_{D;r}$  is the  H\"older
  space; the
 $\Lambda_r(\ov D)$ with norm $|\cdot|_{\Lambda_r(\ov D)}$  is the
  Lipschitz  space  (see~\rd{def3.1}).
  We mention two main features in $H_q$: the first is a  {\it regularized} Henkin-Ram{\'{\i}}rez map,  introduced in this paper, for a strictly pseudoconvex domain with $C^2$ boundary, and the second
  is  the commutator $[\db,E]$,   defined by  $[\db,E]\var=\db E\var-E\db\var$. The commutator has an important property:
$$
[\db,E]f=0, \quad \text{in $\ov D$}.
$$
Combining with the property $[\db, E]\colon \Lambda_{r}(\ov D)\to \Lambda_{r-1}(\cc^n)$,
the commutator is     a {\it smooth} cut-off operator losing one derivative.
We mention closely related previous work.  Lieb-Range~\ci{LR80} first introduced the Seeley extension operator for their $\db$ solution operators in strictly pseudoconvex domains  and the extension has been a basic technique in   other situations.    Ma-Michel~\ci{MM93} used it for a suitable domain in a strictly pseudoconvex real hypersurface in $\cc^n$ for $n\geq4$ and Alexandre~\ci{Al06} used it for finite type convex domains.
If $\var$ is $\dbar$-closed, we obtain  $[\db,E]\var=\db E\var$ and the $H_q$ is an analogue of  the Lieb-Range $\dbar$ solution operator. We should mention that the important commutator $[\db,E]$ was introduced by Peters~\ci{Pe91} and it has been used by Michel~\ci{Mi91}, Michel-Shaw~\ci{MS99}, and others.

 A detailed version of \rt{thm1} (\rt{full})  yields  the following.
\begin{cor}\label{dbsol}  Let $r>1$ and $0<q<n$.
 Let $D$ be a bounded strictly pseudoconvex domain  of $C^{2}$ boundary in $\cc^n$. Let $\var\in\Lambda_r(\ov D)$ be a $\db$-closed $(0,q)$-form in $D$. Then there is a solution $u\in \Lambda_{r+1/2}(\ov D)$ to $\db u=\var$ in $\ov D$.
\end{cor}
   Our $C^{k+1/2}$ estimate   improves the regularity results of Siu~\ci{Si74} (for $q=1$) and  Lieb-Range~\ci{LR80} (for all $q$) for the case when    $k$ is an integer bigger than $1$ and $\pd D\in C^{k+2}$. When $\pd D\in C^\infty$ additionally, the improvement for all $r>0$ was obtained by   Greiner-Stein~\ci{GS77}*{Thm.~16.7(c), p.~174} and Phong-Stein~\ci{PS77}  (for $q=1$), and by  Chang~\ci{Ch89}*{Thm. 4.10 $(iii)$ with $U=D$, $q\geq1$}.
The case $q=n$, which is not included in the corollary, is simple: We  need the domain to be Lipschitz, but not necessarily pseudoconvex, while  solutions gain  a full derivative; see \rp{qisn} for details.


 We will also obtain a   Bochner-Martinelli-Leray-Koppelman formula:  If $f$ is a $C^1$ function in $\ov D$ with   $\db f\in C^1(\ov D)$,  then
$$ 
f=H_0f+H_1\db f,
$$ 
where $D$ is strictly pseudoconvex with $C^2$ boundary  and
$$ 
H_0f = \int_{\cc^n\setminus D}\Om_{0,0}^1 \wedge [\db,E] f.
$$ 
Here $\Om_{0,0}^1$ is a Cauchy-Fantappi\`e form of the above-mentioned regularized Henkin-Ramirez function. In connection with previous work,   $H_0f$ is a   holomorphic projection analogous to
$
\widetilde H_0f = \int_{\pd D}\Om_{0,0}^1 f.
$
We will show  in
  \rt{full}
  that  when  $r>1$ the holomorphic  projection $H_0$ maps $\Lambda_r(\ov D)$ continuously  into itself. For $\widetilde H_0$,  Elgueta~\ci{El80} obtained a similar estimate  with a minor loss of regularity and Ahern-Schneider~\ci{AS79}  obtained a sharp estimate
  that actually holds for all $r>0$.
  See also Phong-Stein~\ci{PS77} for the regularity of Bergman and Szeg\H{o} projections for strictly pseudoconvex domains with $C^\infty$ boundary.

As mentioned earlier, one of  our main results is a homotopy formula in \re{va0qH} and \re{nhq}, which admits H\oo lder estimates in $\ov D$. Using $H_q$, we will study the elliptic differential
$$\cL D:=\dbar_z+d_t$$
in $(z,t)\in \cc^n\times\rr^m$,  introduced by Treves~\ci{Tr92}.  Let $D\times S$ be a product domain in $\cc^n\times\rr^m$.
 A $k$-form $\var$ in $D\times S$ is said of {\it mixed type} $(0,k)$ if
$$
\var(z,t)=[\var]_0(z,t)+\dots+[\var]_k(z,t)
$$
where $[\var]_i$ has type $(0,i)$ in $z$ and degree $k-i$ in $t$. By a  $\cL D$-closed form $\var$, we  mean   $\cL D\var=0$. We have the following.
\th{regD} Let $1\leq r\leq\infty$.  Let $D$ be a bounded strictly pseudoconvex domain with $ C^{2}$ boundary and let $S$ be a bounded star-shaped domain in $\rr^m$. Let $\var$ be a $\cL D$-closed form of mixed type $(0,k)$  in $D\times S$ with $k\geq1$. If   $\var\in C^r(\ov D\times\ov S)$, there is a solution $u\in C^{r}(\ov D\times\ov S)$ satisfying $\cL Du=\var$.
\eth
Hanges and Jacobowitz~\ci{HJ97} proved the interior $C^\infty$ regularity of a $\cL D$-solution on a smooth domain $\Omega$ in $\cc^n\times\rr^m$ under a  strictly   Levi  convex condition.

%
%
%
%
%
%
%

We further mention  some important results concerning $\db$ or $\db_b$ solutions. The $C^\infty$   regularity results of $\db$ solutions were achieved by Kohn~\ci{Ko64} for smoothly bounded strictly pseudoconvex domains, by Kohn~\ci{Ko64} for $n=2$ and Catlin~\ci{Ca87} for pseudoconvex domains of finite D'Angelo type~\ci{Da82}, and by Kohn  for smoothly bounded pseudoconvex domains~\ci{Ko73}. McNeal~\ci{Mc94} obtained exact subelliptic estimates  for finite type convex domains.
The results of finite smoothness solutions have also been obtained. Folland and Stein~\ci{FS76} obtained the regularity in  non-isotropic Lipschitz spaces for $\db_b$ and $\Box_b$ solutions on strictly pseudoconvex CR manifolds.
  The regularity of  $\db$ and $\dbb$ solutions for $(0,1)$-forms was obtained  by Chang-Nagel-Stein~\ci{CNS92}, Fefferman-Kohn~\ci{FK88}, and Christ~\ci{Ch88} for  finite type pseudoconvex domains in $\cc^2$,
  and by Fefferman-Kohn-Machedon~\ci{FKM90} for   finite type domains in $\cc^n$ with diagonalizable Levi-form.
Note that
 Shaw~\ci{Sh91}  obtained the exact $C^{1/m}$ estimate of $\db_b$ solutions for $(0,1)$-forms in the boundary of an ellipsoid of finite type $m$ in $\cc^n$.
   For $(0,q)$-forms, the H\"older estimates for $\db_b$ solutions were finally achieved by Koenig~\ci{Ko02} for finite type CR manifolds with comparable eigenvalues in the Levi form.

\medskip

We now state two   questions.

\medskip

\noindent
{\bf Question 1.} {\em  Let $0< q<n$. Let $D$ be a bounded strictly pseudoconvex domain with $C^2$ boundary
in $\cc^n$ . Let $\var$ be a $\db$-closed $(0,q)$-form in $D$. If $\var\in \Lambda_r(\ov D)$ and $0< r\leq1$, does there exist  $u\in \Lambda_{r+1/2}(\ov D)$  satisfying $\db u=\var$ in $D$?}

\medskip

As mentioned early, when  $\pd D\in C^\infty$,   positive results for the above question via Kohn's solution are in~\cites{GS77, PS77,Ch89} for all  $r>0$. \nrc{dbsol} gives a positive answer when $r>1$.
The result of Kohn~\cites{Ko73}  and \nrc{dbsol} give rise to the following question.

\medskip

\noindent
{\bf Question 2.} {\em  Let $0< q<n$. Let $D$ be a bounded pseudoconvex domain in $\cc^n$ with $C^2$ boundary. Let $\var$ be a $\db$-closed $(0,q)$-form in $D$. If $\var\in C^\infty(\ov D)$, does there exist  $u\in C^{\infty}(\ov D)$  satisfying $\db u=\var$ in $D$?}

\medskip

Finally, we should mention that Chaumat and Chollet~\ci{CC91} obtained a $\db$ solution with a loss of $n-q-1$ derivatives, when $D$ is   a convex domain of $C^2$ boundary and $r\in\nn$. They also obtained other results.   Michel-Shaw~\ci{MS99} also showed that when $D$ is an annulus domain $\Om_1\setminus\ov\Om_2$, where $\Om_1$ is a bounded strictly pseudoconvex domain with   $C^\infty$ boundary, $\Om_2$ is a pseudoconvex domain which is relatively compact  in $\Om_1$ and has $C^2$ boundary, there exists a solution $u\in C^\infty(\ov D)$ to $\dbar u=f$, if $f$ is a  $\db$-closed $(0,q)$-form in $C^\infty(\ov D)$ and   $0< q< n-1$.

The paper is organized   as follows.
In section~\ref{sect:2} we derive the  homotopy formula. In section~\ref{sect:3} we recall the Whitney   and Stein extension operators from \ci{St70} and use them to obtain regularized defining functions for domains with $C^2$ boundary and describe equivalent norms of $\Lambda_r(\ov D)$. Section~\ref{sect:4} contains the main estimation of this paper,  assuming the existence of    regularized Henkin-Ram{\'{\i}}rez functions. The latter are derived in section~\ref{sect:5} for which we follow the classical construction of Henkin-Ram{\'{\i}}rez functions. The final section contains two homotopy formulae for the $\cL D$-complex    and the proof of \rt{regD}.

\medskip

\noindent
{\bf Acknowledgments.} The author is grateful to   Andreas Seeger for helpful discussions on the real  interpolation theory.

  \setcounter{thm}{0}\setcounter{equation}{0}

\section{The homotopy formula and the commutator}
\label{sect:2}

In this section we derive a homotopy formula,  inspired by Lieb-Range~\ci{LR80}, Peters~\ci{Pe91}, and Michel-Shaw~\ci{MS99}. We derive it by keeping the minimum smoothness conditions on the domain and  the forms. We will apply it to prove our main results, after the {\it regularized} Henkin-Ram{\'{\i}}rez functions are constructed in section~\ref{sect:5}.

We  first recall  the  Leray-Koppelman homotopy formula.
Let $D$ be a  bounded domain with $C^{1}$ boundary.  Let $g^1\colon D\times\pd D\to\cc^n$ be a $C^1$   mapping satisfying
$$
g^1(z,\zeta)\cdot(\zeta-z)\neq0, \quad\forall \zeta\in\pd D, z\in D.
$$
  Let
$
g^0(z,\zeta)=\ov\zeta-\ov z
$
and $w=\zeta-z$.      Define  \gan
\omega^i=\f{1}{2\pi i}\f{g^i\cdot dw}{g^i\cdot w},
\quad
\Omega^i=\omega^i\wedge(\ov\pd\omega^i)^{n-1},\\
\Omega^{01}=\omega^0\wedge\omega^1\wedge\sum_{\alpha+\beta=n-2}
(\ov\pd\omega^0)^{\alpha}\wedge(\ov\pd\omega^1)^{\beta}.
\end{gather*}
Here both differentials $d$ and $\db$ are in $z,\zeta$ variables.
 We have
$$ 
\om^i\wedge(\db\om^i)^\all=
\frac{g^i\cdot dw\wedge(\db (g^i\cdot dw))^\all}{(2\pi \sqrt{-1} \, g^i\cdot w)^{\all+1}}, \quad\all=1,2,\dots.
$$ 
We decompose $\Om^i=\sum\Om_{0,q}^i$ and $\Om^{01}=\sum\Om_{0,q}^{01}$,
where $\Om_{0,q}^i$ (resp.~$\Om_{0,q}^{01}$)
  has type $(0,q)$  in $z$ and   type $(n,n-1-q)$ (resp.~$(n,n-2-q)$)  in $\zeta$. Set $\Om_{0,-1}^1=0$ and
   $\Om^{01}_{0,-1}=0$. By the Koppelman lemma~\cite{CS01}*{p.~263}, we have
\ga
\label{kop1}\db_\zeta\Om_{0,q}^1+\db_z\Om_{0,q-1}^1=0,   \quad q\geq 0,\\
\label{kop2}\db_\zeta\Om^{01}_{0,q}+\db_z\Om_{0,q-1}^{01}=\Om_{0,q}^0-\Om_{0,q}^1, \quad q\geq0.
\end{gather}

We  need to know how the sign changes, when the exterior differential interchanges with integration. Following
  notations in Chen-Shaw~\ci{CS01}*{p.~263}, we define
$$
  \int_{y\in M} u(x,y)dy^J\wedge dx^I =  \left \{\int_{y\in M}u (x,y)dy^J\right\}dx^I
$$
for a continuous function $u$ in a manifold $M$. If $d_x$ is the exterior differential in $x$-variables, we have
\eq{checksign}
d_x\int_M\phi(x,y) =(-1)^{\dim M}\int_Md_x\phi(x,y).
\eeq
The  Leray-Koppelman
 homotopy formula \cite{CS01}*{p.~273} for a $(0,q)$-form $\var$ is given by
\ga\label{chtf--}
\var(z)=\db_zT_q\var+T_{q+1}\db_z\var, \quad z\in D, \quad 1\leq q\leq n,\\
\var(z)=\int_{\pd D}\Om_{0,0}^1\var +T_1\db \var,  \quad q=0,
\label{lerayf-}
\intertext{with}
T_q\var=-\int_{\pd D}\Om_{0,q-1}^{01}\wedge \var+\int_D\Om_{0,q-1}^0\wedge \var,
\label{chtf-}
\quad q\geq1.
\end{gather}

  \pr{hf}
Let $D\subset\cc^n$ be a  domain  with $C^{1}$ boundary and let $\cL U$ be a bounded neighborhood of $\ov D$.  Let $g^0(z,\zeta)=\ov\zeta-\ov z$.  Let $g^1(z,\zeta)=W(z,\zeta)$ where $W\in C^{1}(D\times(\cL U\setminus D))$
is a Leray mapping, that is that $W$ is holomorphic in $z\in D$ and satisfies
$$ 
\Phi(z,\zeta):=W(z,\zeta)\cdot(\zeta-z)\neq0, \quad z\in D, \quad \zeta\in\cL U\setminus D.
$$ 
Let $\var$ be  a $(0,q)$-form in $\ov D$. Suppose that $ \var$ and $\db\var$ are in $C^1(\ov D)$. Then in $D$
\ga\label{tsqf}
\var=\db H_q\var+H_{q+1}\db\var, \quad 1\leq q\leq n,\\
\label{tsqf+}
\var=H_0\var+H_1\db \var, \quad q=0,\\
\intertext{where}
\label{Hqv}
H_q\var:=\int_{\cL U}\Om_{0,q-1}^0\wedge E\var+\int_{\cL U\setminus D}\Om_{0,q-1}^{01}\wedge[\db,E]\var,  \quad q>0,\\
H_0\var:=\int_{\pd D}\Om_{0,0}^1\var-\int_{\cL U\setminus D}\Om_{0,0}^1\wedge E\db \var=\int_{\cL U\setminus D}\Om_{0,0}^1\wedge [\db, E] \var.
 \label{H0f}
\end{gather}
\epr
\begin{proof} In the formulae,  the extension $E$ constructed in ~\ci{St70} will be recalled in \rl{steinext} below.
The   $E$ is defined for functions. We thus define $E\var$ by applying $E$ componentwise to its coefficients, which results in a form of the same type. We always assume that $E\var $ has a compact support in $\cL U$, by using a cut-off function.

 Assume that $q\geq1$. Let us modify the solution operator $T_q$   given by \re{chtf--}-\re{chtf-}, by applying the method of Lieb-Range~\ci{LR80} via  the linear extension $E$.
The $\Om^{01}$ has total degree $2n-2$. Applying Stokes' formula and \re{kop2}-\re{checksign}, we   get
\al\label{keyid}
&-\int_{\zeta\in\pd D}\Om_{0,q-1}^{01} \wedge \var =\int_{\zeta\in\cL U\setminus D}\Om_{0,q-1}^{01}\wedge\db_\zeta E\var
+\int_{\zeta\in\cL U\setminus D}\db_\zeta\Om_{0,q-1}^{01}\wedge  E\var
\\
&\qquad\quad  =\int_{\cL U\setminus D}\Om_{0,q-1}^{01}\wedge\db E\var
\nonumber \\
\nonumber&\qquad \qquad -\int_{\cL U\setminus D}\left(\db_z\Om_{0,q-2}^{01}\wedge E\var+\Om_{0,q-1}^1\wedge E\var  -\Om_{0,q-1}^0\wedge E\var \right)
\\
\nonumber
&\qquad \quad =\int_{\cL U\setminus D}
\Om_{0,q-1}^{01}\wedge\db E\var
-\db_z\int_{\cL U\setminus D}\Om_{0,q-2}^{01}\wedge E\var \\
\nonumber
&\qquad \qquad+\int_{\cL U\setminus D}
\left(- \Om_{0,q-1}^1\wedge E\var+
\Om_{0,q-1}^0\wedge E\var \right).
\end{align}
Let us apply $\db$  to the last $4$ terms. The  second of the four terms becomes zero. The third also becomes zero since it is   holomorphic   for $q=1$ and it is zero for $q>1$. Thus we obtain for $z\in D$
\al\label{keyidsim}
&-\db\int_{\zeta\in\pd D}\Om_{0,q-1}^{01}( z,\zeta )\wedge \var(\zeta)+\db \int_{\zeta\in D}\Om_{0,q-1}^0( z,\zeta )\wedge \var(\zeta)
 \\
 &
 \nonumber \qquad\qquad=\db\int_{\cL U \setminus D}
\Om_{0,q-1}^{01}(z,\zeta)\wedge\db E\var(\zeta) +
\db\int_{\cL U  }
\Om_{0,q-1}^0(z,\zeta)\wedge E\var(\zeta).
\end{align}
So far, we have used $\var\in C^1(\ov D)$. Assume now that $\db\var\in C^1(\ov D)$. Using the last 4 terms in \re{keyid} in which $\var$ is replaced by $\db\var$, we obtain
\al\label{keyidsim+}
-\int_{\pd D}\Om_{0,q}^{01}\wedge\db \var &+\int_D\Om_{0,q}^0\wedge \db \var
=\int_{\cL U\setminus D}
\Om_{0,q}^{01}\wedge\db E\db \var
\\
&- \db
\int_{  {\cL U}  \setminus D}
\Om_{0,q-1}^{01}\wedge E\db \var-\int_{  {\cL U}  \setminus D}
 \Om_{0,q}^1\wedge E\db \var \nonumber\\
 & +\int_{  {\cL U}  \setminus D}
\Om_{0,q}^0\wedge E\db\var +\int_{D }
\Om_{0,q}^0\wedge \db \var.
\nonumber  \end{align}
On the right-hand side, the first term can be written via the commutator as
 $\db E\db \var =(\db E-E\db)\db \var$.
 Since $q\geq1$,  the third  is zero. The second, when combined with the first  term  on the right-hand side of \re{keyidsim},
 gives us the desired commutator for $\var$.    Adding   \re{keyidsim}-\re{keyidsim+} yields \re{tsqf}.

 To derive \re{H0f},  we recall that  $\Omega_{0,-1}^{01}=0$ and by \re{keyidsim+} we get
 \aln
-\int_{\pd D}\Om_{0,0}^{01}\wedge\db \var &+\int_D\Om_{0,0}^0\wedge \db \var
=\int_{\cL U\setminus D}
\Om_{0,0}^{01}\wedge\db E\db \var
\\
& -\int_{  \cL U  \setminus D}
 \Om_{0,0}^1\wedge E\db \var  +\int_{  \cL U   }
\Om_{0,0}^0\wedge E\db \var=H_1\db\var-\int_{  \cL U  \setminus D}
 \Om_{0,0}^1\wedge E\db \var.
\nonumber  \end{align*}
Thus we have verified \re{tsqf+} with $H_0$ being defined in \re{H0f}, while the second expression of $H_0$ in \re{H0f} follows from Stokes' formula and
$\db_\zeta\Om_{0,0}^1=0$ by \re{kop1}.
\end{proof}

Throughout the paper, $|\cdot|_{\ov D;r}$, or $|\cdot |_r$ for abbreviation,  denotes the  H\"older $C^r$ norm, $r\in[0,\infty)$,  for differential forms or functions on a domain $D$.
 We finish the section with the following  interior  estimate of Webster~\ci{We89}.
\pr{estL} Let $r\in[0,\infty)$. Let $\cL U$ be a bounded domain in $\cc^n$ with $\ov D\subset\cL U$.  Let $L\psi=\int_{\cL U}\Om_{0,q}^0\wedge \psi$ with $0\leq q\leq n$. Then
\ga
|L\psi|_{\ov D;r}\leq C^*_{a}|\psi|_{\cL U;r},\\
\label{Wea+1}
|L\psi|_{\ov D;r+1}\leq C^*_r  |\psi |_{\cL U;r}, \quad r\not\in\nn,
 \end{gather}
where $C^*_r \leq C_r(\cL U)\dist(D,\pd \cL U)^{-c_0r-c_1}$ and $C_r(\cL U)$ depends only on $r$ and the diameter of $\cL U$.
\epr


\setcounter{thm}{0}\setcounter{equation}{0}

\section{Regularized  defining functions and preliminaries for Lipschitz estimates.}
\label{sect:3}

In this section  we define  a  {\it regularized} defining function  for a domain $D$ by Whitney's extension so that the derivatives of the extension  have optimal growth rates near the boundary of $D$. The defining function will play an important role in our estimates. We  recall an extension operator of Stein~\ci{St70} and basic facts about the Lipschitz space $\Lambda_r$ and its equivalent norms. The equivalent norms   are used for the $\Lambda_{r+1/2}$ estimate when $r+1/2=2,3,\dots$.
 We will also recall some basic results on the real interpolation theory. The interpolation will be used for $C^{r+1/2}$ estimates when $r=2,3,\dots$, which as mentioned in the introduction,
  improves the regularity result of Lieb-Range.
 While the results of this sections might be known to the reader, we formulate them for the purpose of this paper. We will also specify the dependence of the  various constants  on the domains, which is used to address the stability of estimates of the homotopy operators in \rt{full}.
 We will conclude the section with a regularity result for the $\db$ equation of top type.

Let us first introduce   notations. For $r\in\rr$, $[r]$ denotes the largest integer $k\leq r$.
For two sets $A,B$ in $\rr^n$, $\dist(A,B)$ denotes $\inf\{|a-b|\colon a\in A,b\in B\}$.
For $\all\in \nn^n$, let $$\pd_x^\all f:=\pd_{x_1}^{\all_1}\dots\pd_{x_n}^{\all_n}f(x)$$
 denote the partial derivative function in $x$
  and $\pd_x^kf$ also denotes the set of all partial derivatives of order $k$.
Let $D$ be a domain in $\rr^n$.    Let $C^r(\ov D)$ denote the set of functions $f$ in $D$ such that $\pd_x^\all f$ extend to   functions $f^{(\all)}\in C^{r-[r]}(\ov D)$ for  $|\all|:=\all_1+\dots+\all_n\leq{r}$. For a continuous function $f$ in $\ov D$, $|f|_{\ov D;0}$ denotes $|f|_{L^\infty(D)}$, while $|f|_{\ov D;r}$ also denotes the H\"older norm for $f\in C^r(\ov D)$ when $0<r<1$. For $f\in C^r(\ov D)$, define
$$
|f|_{\ov D;r}:=\max_{|\all|\leq r} |f^{(\all)}|_{\ov D;r-[r]}.
$$

\subsection{Regularized defining functions.}
Let $F$ be a closed set in $\rr^n$ and let $r\in(0,\infty)$.   We recall the following definition. 
\begin{defn}(\cite{Wh34}*{p.~64},\ci{St70}*{p.~194})\label{WhitneyCa} Let $F$ be a non-empty closed subset of $\rr^n$.
A function $f$ in $F$ is said in $C_w^r(F)$ in terms of the functions $f^{(\all )}$ in $F$ for $\all \in\nn^n$ and $|\all |\leq{r}$,  if $f^{(0)}=f$ and  there is a finite constant $A$ so that $
|f^{(\all )}(x)|\leq A$ for all $x\in F$, while $R_\all$, defined by
\eq{defRj}
P_\all (x,p):=\sum_{|\beta |+|\all |\leq{r}}\f{f^{(\all +\beta )}{(p)}}{\beta !}(x-p)^\beta , \quad R_\all (x,p):=f^{(\all )}(x)-P_\all (x,p)
\eeq
has the properties: (i) $|R_\all (x,p)|\leq A|x-p|^{r-|\all |}$ for $x,p\in F$ and $|\all|\leq r$; (ii)  when $r\in\nn$,
 for each $p\in F$ and $\e>0$ there is $\del>0$ so that for $x,x'\in F$ with $|x-p|+|x'-p|<\del$,
\eq{eRj} |R_\all (x',x)|\leq \e|x'-x|^{[r]-|\all |}.
\eeq
As observed by Whitney~\ci{Wh34}, condition \re{eRj} is essential and consequently all $f^{(\all)}$ are continuous in $F$
for $|\all|\leq r$.
Following Stein~\ci{St70}*{p.~173},
we define $|f|^w_{F;r}$ to be the  infimum of the constants $A$ for all possible choices of $f^{(\all )}$ for $0\leq|\all |\leq{r}$.
\end{defn}

\pr{Whit} Fix $r\in[0,\infty)$. Let $F$ be a closed subset of $\rr^n$.    There is an extension operator $E_r \colon C_w^r(F)\to C^r(\rr^n)$    so that $|E_r  f|_{\rr^n;r}\leq C_r|f|^w_{F;r}$. Moreover,    for $x\in F^c:=\rr^n\setminus F$ with $d(x):=\dist(x,F)<1$,
\ga
\label{pdJe}
|\pd_x^\all E_rf-P_{\all}(x,x_*)|\leq C_r\sup_{|\beta|\leq r,x'\in\pd F, |x'-x_*|<4d(x)}|R_{\beta}(x^*,x')| d(x)^{|\beta|-|\all |},\quad|\all|\leq r;\\
|\pd_x^k E_rf|\leq C_{r} |f|^w_{\ov D;r}(1+d(x)^{r-k}),\quad x\in F^c, \quad k=0,1,2,\dots,\label{nabE}
\end{gather}
where  $|x-x_*|=d(x)$ with $x_*\in F$, and $P_\beta=0$ for $|\beta|>r$.
\epr
\begin{proof} Inequality \re{nabE} is proved in Stein~\cite{St70}*{p.~178} and Glaeser~\ci{Gl58} when $|\all|\leq{r}+1$ and stated in~\ci{Gl58}*{p.~31} for all $\all$. We present here a proof for the reader's convenience.  Recall  from   \cite{St70}*{p.~169-170}  the following properties: $(i)$~There are $\var_k\in C^\infty_0(F^c)$ so that $\sum\var_k=1$ in $F^c$,  $0\leq\var_k\leq1$, and
\ga
|\pd_x^\all\var_k|\leq C_\all\dist(x,F)^{-|\all|},\quad \supp\var_k\subset Q_k,
\end{gather}
where $Q_k$ are cubes satisfying $\f{1}{2}\diam Q_k\leq\dist(Q_k,F)\leq 5\diam Q_k$, and $F^c=\cup_kQ_k$.
$(ii)$ Each point in $F^c$ is contained in at most $N_0$ of cubes $Q_k$. Here $N_0$ and $C_\all$ are independent of $F$.
%

For each $Q_k$, fix $p_k\in F$ such that $\dist( F,Q_k)=\dist(p_k,Q_k)$. We choose $\{f^{(\all)}\colon|\all|\leq r\}$ so that the constant $A$ in \rda{WhitneyCa} satisfies $A\leq 2|f|^w_{r}$. Let 
$P(x,p)=\sum_{|\all|\leq r}\f{1}{\all!}f^{(\all)}(p)(x-p)^\all$. Define $E_{r}f=f$ in $F$ and
\gan
E_rf(x)=\sum_i\nolimits'P(x,p_i)\var_i(x), \quad x \in F^c.
\end{gather*}
Here the sum with the prime is   over the $i$ satisfying $\dist (Q_i,F)<1$. When $d(x)<1$ and $x\in Q_i$, we have $\dist(Q_i,F)<1$.  Thus we   drop the prime, by assuming   $d(x)<1$.  Recall from \cite{St70}*{Lemma, p.~177} that
\ga  P_\beta(x,p)-P_\beta(x,q)
=\sum_{|\gamma|\leq r-|\beta|}R_{\beta+\gamma}(p,q)\f{(x-p)^\gamma}{\gamma!},
\quad p,q\in F,\label{PaPb}\\
 \pd_x^\beta P(x,p)=P_\beta(x,p), \quad p\in F.
\end{gather}
For $x\in F^c$, we fix $x_*\in F$ such that $|x-x_*|=\dist(x,F)$.  Suppose that $|\all|>0$.  Then
\al\label{JEf}
\pd_x^\all E_rf&=\sum_{\beta,\all-\beta\in\nn^n}\sum_k\binom{\all}{\beta}\pd^\beta P(x,p_k)\pd_x^{\all-\beta}\var_k\\
&=P_\all(x,x_*)+
\sum_{\beta,\all-\beta\in\nn^n}\sum_k\binom{\all}{\beta}\left(P_\beta(x,p_k)-P_\beta(x,x_*)\right)\pd_x^{\all-\beta}\var_k.\nonumber
\end{align}
We  only need to consider the terms with $\pd_x^{\all-\beta}\var_k\neq0$. Thus $x\in Q_k$ and   there are  at most $N_0$ of such $\var_k$'s. When $ \pd_x^{\all-\beta}\var_k\neq0$, we have
 $|x-p_k|\leq \diam(Q_k)+\dist(Q_k,F)\leq 3\dist(Q_k,F)\leq 3d(x)$.  Then by \re{PaPb}
\al\label{estadd}
&|P_\beta(x,x_*)-P_\beta(x,p_k)||\pd_x^{\all-\beta}\var_k|
  \leq C \sum_{|\gaa|\leq r-|\beta|}|R_{\beta+\gaa}(x_*,p_k)| d(x)^{|\beta|+|\gaa|-|\all|}.
\end{align}
Combining it with \re{JEf} yields
$$
|\pd_x^\all Ef-P_{\all }(x,x_*)|\leq C \sum_{0\leq|\beta|\leq r}  |R_{\beta}(x_*,p_k)| d(x)^{|\beta|-|\all|}.
$$
We also have $|x_*-p_k|\leq|x_*-x|+|x-p_k|\leq 4d(x).$
Hence $|\pd_x^\all Ef-P_{\all}(x,x_*)|\leq C' |f|^w_r d(x)^{r-|\all|}$. Also the continuity of $\pd^\all E_af$ comes from
$$
|\pd_x^\all Ef-P_{\all}(x,x_*)|\leq C'\e_{x,x*}d(x)^{[r]-|\all|}
$$
with $\e_{x,x^*}\to 0$ as $x$ tends to $x_0\in \pd D$.
We have proved \re{pdJe}, while \re{nabE} follows directly from \re{pdJe}.  When $r\in\nn$, \re{nabE} implies that $|E_rf|_r\leq C|f|_{r}^w$.
When $r>[r]$, \re{nabE}
for $k=[r]+1$ also implies that $|\pd ^{[r]}E_rf |_{\rr^n;r-[r]}\leq C|f|^w_{r}$; see the proof in~\ci{St70}*{Thm.~3, p.~173}.
\end{proof}

\begin{rem}When $F=\ov D$ for a domain $D$ in $\rr^n$ and $f\in C_w^r(F)$, $D$ is dense in $F$ and   $f^{(\all)}$ are uniquely determined by the values of $f$ in $D$; in fact $f^{(\all)}=\pd^\all f$ in $D$. In this case, the above $E_r $ is a linear operator for a fixed sequence $p_k$ appeared in the above proof.
\end{rem}

We first identify $C^r(\ov D)$ with $C^r_w(\ov D)$ under a mild condition on the domain $D$.
\le{estR}Let $r\geq1$ and $L\geq1$.  Let $D$ be a domain in $\rr^n$. Assume   that any two points $p,q$ in $\ov D$ can be connected by $\gaa$, a union of finitely many line segments in $\ov D$, so that $\gaa$ has length at most $ L|p-q|$ and $\gaa\cap\pd D$ is a finite set.
Then $C^r(\ov D)=C^r_w(\ov D)$   and
 $|f|_{\ov D;r}\leq |f|_{\ov D;r}^w\leq C_rL^r|f|_{\ov D;r}$.
\ele
\begin{proof} When $|\all |=[r]$, we have $R_\all (x,y)=\pd^\all _xf-\pd^\all _yf$. By the continuity of $\pd_x^\all f$ and the definition of H\oo lder ratio, we get $|R_\all(x,y)|\leq |f|_{D;r}|x-y|^{r-[r]}$ and  get \re{eRj} by the continuity of $\pd^\all f$.

Assume that $|\all |<[r]$.
Let $\gaa\colon[0,1]\to\ov D$ be a piecewise linear curve with $\gaa(0)=p,\gaa(1)=q$. Suppose that $\gaa(t)\in D$ and $|\gaa'(t)|\leq L|p-q|$ for $t\in(t_k,t_{k+1})$ with $t_0=0,\dots, t_N=1$. Choose an
increasing $C^\infty$ function $\hat s$ such that $\hat s(0)=0$,  $\hat s(1)=1$, and all derivatives of $\hat s$ vanish at $0,1$.  Let $s(t)=t_k+(t_{k+1}-t_k)\hat s((t-t_k)/{(t_{k+1}-t_k)})$ for $t\in[t_k,t_{k+1}]$. Then $s(t_j)=t_j$, $s^{(\ell)}(t_j)=0$ for all $\ell>0$, and $0\leq s'(t)\leq C$, where $C$ is independent of $t_i, N$. Then $t\mapsto\gaa(s(t))$ is a $C^\infty$ curve connecting $p,q$. Let $\gaa(t)$ still denote $\gaa(s(t))$. We have $|\gaa'(t)|\leq CL|p-q|$.

Let $g(t)=R_\all (\gaa(t),p)$.   Then $g$ is in $C^1([0,1])$.  We have
$$
g(1)=\sum_i\int_0^1\pd_{x_i}|_{x=\gaa(s_1)}R_\all (x,p)  \gaa_i'(s_1)\, ds_1.
$$
Also,  $\pd_x^\beta R_\all (x,p)$ vanishes at $x=p$ if $|\beta|+|\all |\leq r$. Therefore, 
$$
g(1)=\sum \int_0^1\cdots\int_0^{s_{i-1}} \pd_{x_{k_i}}\cdots\pd_{x_{k_1}}R_\all (\gaa(s_i),p)\gaa_{k_1}'(s_1)\cdots\gaa_{k_i}'(s_i)\, ds_i\dots ds_1,
$$
for summing over $k_1,\dots, k_i$ with  $k_1+\dots+k_i=i$ and $i=[r]-|\all |$.
We obtain $$|R_\all (q,p)|\leq C_i(CL|p-q|)^{[r]-|\all |}\max_{t\in[0,1],|\beta|=[r]-|\all |}|\pd_x^\beta R_\all (\gaa(t),p)|.
$$
Note that $\pd_x^\beta R_\all (x,p)=R_{\beta+\all } (x,p)$.  The lemma is verified.
\end{proof}

The proof also yields the following inequality.
\pr{remainder}Let $D$ be as in \rla{estR}. Let $P(x,p)$ be the Taylor polynomial of $f$ of degree $k$ about $p\in\ov D$.  Then for $x\in\ov D$
$$
|f(x)-P(x,p)|\leq C_k L^k|x-p|^{k}\sup_{x',|\all|=k} |\pd^{\all}f(x')-\pd^{\all}f(p)|,
$$
where   $x'\in \ov D$ and $|x'-p|\leq L|x-p|$.
\epr

\begin{defn}[{\cite{St70}, p.~189}]\label{mins}
Let $D$ be a domain in $\rr^n$. We say that $\pd D$ is    {\it minimally smooth} if the following conditions hold: There are positive numbers $\e,N,M$, and a sequence of open subsets $U_1,U_2,\dots$ of $\rr^n$  so that the following hold:
\bppp
\item If $x\in\pd D$, then $B(x,\e)\subset U_i$; $B(x,\e)$ is the ball of center $x$ and radius $\e$.
    \item No point of $\rr^n$ is contained in more than $N$ of the $U_i's$.
    \item For each $i$ there exists a domain $D_i$ in $\rr^n$, defined by
    $
    x_{n_i}>\var_i(x_{n_i}')
    $
 for $x=(x_1,\dots, x_n)$,  $x_{n_i}'=(x_1,\dots,\hat x_{n_i},\dots x_n)$
 so that $U_i\cap D=U_i\cap D_i$ and
$$
|\var_i(u)-\var_i(v)|\leq M|u-x|, \quad u,v\in\rr^{n-1}.
$$
\eppp
We will denote by $C_r(D)$ a finite number depending on the above $M,N,\e$, and $r$.\end{defn}
%
%

Note that a bounded domain in $\rr^n$ has a (strong) Lipschitz boundary, i.e. its boundary is locally the graph of a Lipschitz function in some smooth coordinates,  if and only if its boundary is minimally smooth.
\le{regdef}  Let $D$ be a bounded domain in $\rr^n$ with $C^2$ boundary. Let $\rho_0\in C^2(\ov D)$ with $\pd \rho_0\neq0$ in $\pd D$ and $\rho_0\leq0$ in $\ov D$.  There exists a real function
$\rho\in C^2(\rr^n)\cap C^\infty(\rr^n\setminus\ov D)$ such that $\rho=\rho_0$  in $\ov D$, and for $0<d(x):=\dist(x,D)<1$,
\ga \label{regca}
|\pd_x^i\rho|\leq C_{i} L^2 |\rho_0|_{\ov D;2}(1+d(x)^{2-i}),\quad   i=0,1,2,\dots, \\
\label{pdrpdr}
|\pd \rho(x)-\pd \rho_0(x_*)|\leq CL^2|x-x_*|\max_{y\in\ov D,|y-x_*|\leq4L|x-x_*|}|\pd_y^2\rho_0|, \\
|\pd^2\rho(x)-\pd^2\rho_0(x_*)|\leq CL^2\om_2(|x-x_*|),
\label{pd2r}
\end{gather}
where  $x_*\in \pd D$,   $|x-x_*|= \dist(x,D)$, and
\gan
\om_2(\del)=\sup_{x'\in\ov D,x\in\pd D,|x'-x|\leq 4L\del}|\pd^2\rho_0(x')-\pd^2\rho_0(x)|.
\end{gather*}
If $0<d(x)<\min_{y\in\pd D}\{1, |\pd_y\rho_0|/{(CL^2|\pd^2\rho_0|_{0})}\}$, then
\eq{rxge}
|\pd_x\rho|\geq \f{1}{2} \min_{y\in\pd D}|\pd_y \rho_0|,\quad
\rho(x)\geq \f{1}{2}\min_{y\in\pd D}|\pd_y \rho_0| d(x).
\eeq
\ele
\begin{proof} Applying \re{pdJe} and \rp{Whit} to $\rho=E_2\rho_0$ and $F=\ov D$,  we obtain
\gan
|\rho(x)|\leq C|(\pd\rho,\pd^2\rho)(x_*)||x-x_*|+C\sup_{x',|\all|\leq 2}|R_\all(x',x_*)|d(x)^{|\all|},\\
|\pd\rho(x)-\pd\rho(x_*)|\leq C|\pd^2\rho(x_*)||x-x_*|+C\sup_{x',|\all|\leq 2}|R_\all(x',x_*)|d(x)^{|\all|-1},\\
|\pd^2\rho(x)-\pd^2\rho(x_*)|\leq C\sup_{x',|\all|\leq 2}|R_\all(x',x_*)|d(x)^{|\all|-2},
\end{gather*}
where $x'\in\pd D$ and $|x'-x_*|\leq4|x-x_*|$. Here $R_\all(x',x_*)$ is the Taylor remainder of $\rho$ defined by \re{defRj} with $r=2$. By \rp{remainder}, we have for $|\all|\leq2$
$$|R_\all(x',x_*)|\leq C(L|x'-x_*|)^{2-|\all|}\sup_{x''\in\ov D,|x''-x_*|\leq L|x'-x_*|}|\pd^2\rho(x'')-\pd^2\rho(x_*)|.$$
This gives us \re{pdrpdr}-\re{pd2r}. We get \re{regca}  from \re{nabE} and $|E_2r_0|_2\leq CL^2|r_0|_{\ov D;2}$ by  \rl{estR}.  Estimate \re{rxge} follows directly from \re{pdrpdr}.
\end{proof}
We   call the above $\rho$ a   {\it regularized}  $C^2$ defining function of $D$.
We will also need the following version of Stokes' theorem.
\le{regstokes} Let $m$ be a positive integer, and let $b\in\rr$.  Let $\cL V $ be a bounded domain in $\rr^n$ with $C^1$ boundary. Assume  that $B$ and $S$ are functions in $   C^1(\cL V)$ and for $x\in\cL V$ and $i=0,1$,
 $$
  |\pd_x^iS|<C\dist(x,\pd \cL V)^{m-i}, \quad
|\pd_x^iB|\leq C(1+\dist^{b-i}(x,\pd\cL V  )),\quad C=C(B,S)<\infty.$$
Assume further that $b+m>0$. Then
$
\int_{\cL V }  B(x)\pd_{x_j} S\, dx=-\int_{\cL V }  S(x)\pd_{x_j} B\, dx.
$
\ele
\begin{proof} Let $d(x)=\dist(x,\pd \cL V)$. We know that   $B \pd_{x_j} S$ is   Lebesgue integrable in $\cL V $, since
$| B(x)\pd_{x_j}S|\leq C'(1+\dist(x)^{m+b-1}).
$
Analogously, $S(x)\pd_{x_j}B$ is integrable in $\cL V $.

Let $N_{\del}$ be the set of $x\in\cL V$ satisfying $d(x)<\del$.
 Take $\chi_\ell\in C_0^\infty(\cL V\setminus N_{1/\ell})$ so that $0\leq\chi_\ell\leq1$, $\chi_\ell=1$ in $\cL V\setminus N_{2/\ell}$,  and $|\pd_{x_j}\chi_\ell|\leq C d(x)^{-1}$.  Then we have
$$
\left|\int_{\cL   V }\pd_{x_j} ((1-\chi_\ell )BS)\, dx\right|\leq C  \int_{N_{2/\ell}} (1+d(x)^{m-1+b})\, dx,$$
which tends to $0$ as $\ell\to\infty$ as $\cL V$ is bounded and $\pd \cL V\in C^1$. Also,
$ \int_{ \cL V } ((1-\chi_\ell) S(x)\pd_{x_j}B\, dx $,  $  \int_{ \cL V }((1-\chi_\ell )B)(x) \pd_{x_j}S\, dx $,
and $\int_{\cL V} (SB)(x)\pd_{x_j}\chi_\ell\, dx$ tend to $0$ as $\ell\to\infty$. By Stokes' theorem, we have
$$
\int_{\cL V} ( SB)(x)\pd_{x_j}\chi_\ell\, dx+\int_{\cL V} ( \chi_\ell B)(x)\pd_{x_j}S\, dx+
\int_{\cL V} ( \chi_\ell S)(x)\pd_{x_j}B\, dx=0.
$$
Letting $\ell\to\infty$, we get the identity.
 \end{proof}

\subsection{Equivalent norms.}
 For $0< r\leq 1$,  the  Lipschitz space $\Lambda_r(\rr^n)$ is the set of functions $f\in L^\infty(\rr^n)$ such that  \gan
|f|_{\Lambda_r(\rr^n)}:=|f|_{\rr^n,r}:=|f|_{L^\infty(\rr^n)}+\sup_{y\in\rr^n}\f{|\Delta_yf|_{L^\infty(\rr^n)}}{|y|^r},\quad 0<r<1,\\
|f|_{\Lambda_1(\rr^n)}:=|f|_{L^\infty(\rr^n)}+\sup_{y\in\rr^n\setminus\{0\}}
\f{|\Delta^2_yf|_{L^\infty(\rr^n)}}{|y|}.
\end{gather*}
Here $\Del_yf(x):=f(x+y)-f(x)$ and thus $\Del^2_yf(x)=f(x+2y)+f(x)-2f(x+y)$.
When $r>1$, we define $\Lambda_r(\rr^n)$ to be the set of functions $f\in C^{[r]-1}(\rr^n)$ satisfying
$$
|f|_{\Lambda_r(\rr^n)}:=|f|_{L^\infty(\rr^n)}+|\pd f|_{\Lambda_{r-1}(\rr^n)}<\infty.
$$
By~\ci{St70}*{Prop.~8, p.~146}, the $|f|_{\Lambda_r(\rr^n)}$ is equivalent with the expression
$$
|f|_{L^\infty(\rr^n)}+\sup_{y\in\rr^n}\f{|\Delta^2_yf |_{L^\infty(\rr^n)}}{|y|^r},
$$
for $0<r<2$.  For a non-integer $r$,   $|\cdot|_{\Lambda_r}$ is equivalent with the H\"older norm $|\cdot|_{\rr^n;r}$ by~\ci{St70}*{Prop.~8, p.~146}.

\begin{defn}\label{def3.1} Let $F$ be a closed subset in $\rr^n$.  Let $r\in(0,\infty)$.
We write $f\in\Lambda_r(F)$ if there exists   $\tilde f\in\Lambda_r(\rr^n)$ such that $\tilde f|_F=f$. Define $|f|_{\Lambda_r(F)}$  to be the infimum of  $|\tilde f|_{\Lambda_r(\rr^n)}$ for all such extensions $\tilde f$.
Sometimes $|f|_{\Lambda_r}$ denotes $|f|_{\Lambda_r(F)}$   for abbreviation.
\end{defn}

We now discuss equivalent norms of the spaces $\Lambda_r$.
The following lemma is   in  McNeal-Stein~\ci{MS94}. We need a quantified version.

\pr{gs141} Let $0<r<\infty$.   Then $f\in\Lambda_r(\rr^n)$ if and only if there is a decomposition
$
f=\sum_{k\geq0} f_k
$
so that $f_k\in C^\infty(\rr^n)$ and
\eq{nabi}
 |\pd^i f_k|_{L^\infty(\rr^n)}\leq A2^{k(i-r)}, \quad i=0,\dots, [r]+1.
\eeq
Furthermore,  the smallest constant $A_r (f)$ of $A$ is equivalent with
$|f|_{\Lambda_r}$,
 i.e. $c_rA_r (f)\leq|f|_{\Lambda_r}\leq C_rA_r (f)$ for some positive numbers $c_r ,C_r$ independent of $f$.
\epr
\begin{proof}The lemma   is proved by Greiner-Stein~\ci{GS77}*{p.~142} for $0<r<1$. For $0< r\leq1$, the existence of decomposition is proved in~\ci{GS77}*{p.~145}. The decomposition is also valid for $r>1$ since $f_k(x)= \int\var_k(t)f(x-t)\, dt$,  each $\var_k\in C^\infty(\rr^n)$ is rapidly decreasing, and hence $\pd_{x_j}f_k(x)=\int\var_k(t)\pd_{x_j}f(x-t)\, dt$.

 Assume that $1\leq r<2$ and \re{nabi} holds.
  We  have $|f|_{L^\infty}\leq \sum A2^{-rk}\leq C_1A$. We decompose
$$
\Del_y^2f(x-y)=\sum_{k\leq N}\Del_y^2f_k(x-y)+\sum_{k>N}\Del_y^2f_k(x-y).
$$
The two sums are bounded by $$
|y|^2\sum_{k\leq N}|\pd^2f_k|_{L^\infty}\leq C_r A |y|^2 2^{2N-rN}, \quad |y| \sum_{k>N}|\pd f_k|_{L^\infty}\leq A|y| 2^{N-Nr}.
$$
When $0<|y|<1$, choose a positive integer $N$ so that $1\leq|y|2^N\leq2$.   Hence $|f|_{\Lambda_1}\leq CA$. Analogously, we can verified the proposition for $r\geq2$.
\end{proof}

We will use a linear extension operator from Stein\ci{St70} to prove the following.
\pr{steinext} Let $D$ be a   domain in $\rr^n$ where $\pd D$ is minimally smooth.
There is a continuous linear extension operator
$E\colon C^0(\ov D)
\to C^0 (\rr^n)$
 so that $Ef=f$ on $D$
 \ga\label{EL}
 |Ef|_{\Lambda_r(\rr^n)}\leq C_r(D)|f|_{\Lambda_r(\ov D)},\quad \forall r\in(0,\infty),
 \\
  \label{EH}
|Ef|_{C^r(\rr^n)}\leq C_r(D)|f|_{C^r(\ov D)},  \quad \forall r\in[0,\infty). \end{gather}
\epr
\begin{proof}We follow the proof in~\ci{St70}. We first recall an extension for each $D_i$. To simplify the notation, we drop the $i$ and assume $n_i=n$. Thus $D$ is defined by $x_n>\var(x')$.  Set $D^c=\rr^n\setminus D$,  ${\ov D}^c=\rr^n\setminus \ov D$, and let $ d(x)$ be the distance of $x$ from $D$. By~\ci{St70}*{Thm.~2, p.~171},  there is a regularized distance function $\Del\in C^\infty(\ov D^c)$    vanishing on $\pd D$ so that
\eq{sDel}
c_1 d(x)\leq\Delta(x)\leq c_2 d(x), \quad |\pd_x^\all\Delta(x)|\leq C_\all\delta^{1-|\all|}(x), \quad x\in \ov D^c.
\eeq
Then we  choose a rapidly decreasing function $\psi\in C^\infty([1,\infty))$ so that
\eq{intp1}
\int_1^\infty\psi(\la)\, d\la=1, \quad\int_1^\infty\la^k\psi(\la)\, d\la=0, \quad k=1,2,\dots.
\eeq
We have a linear extension operator
$$ 
\mathfrak Ef(x)=\int_1^\infty f(x',x_{n}+\la \del^*(x))\psi(\la)\, d\la, \quad x\in \ov D^c,
$$ 
where $\del^*(x)=c\Del(x)\geq 2(\var(x')-x_{n})$.
We need the following estimate:
\eq{eSE}
|\mathfrak Ef|_{L_k^\infty(\rr^n)}\leq C_{k}|f|_{L_k^\infty(D)},
\eeq
where $L_k^\infty(D)$ is the set of functions $f$ in $D$ such that the distributional derivative $\pd^if\in L^\infty(D)$ for all $i=0,\dots, k$, and  the constant $ C_{k}$ depends only on the upper bound of $M$ and $k, p$; see \ci{St70}*{Thm.\! $5^\prime$, p.~181}.

Assume now that $f\in\Lambda_r(\ov D)$. We verify \re{EL} by using an argument in Greiner-Stein~\ci{GS77}*{p.p.~146--147} in which $D$ is a half-space. By the definition of $\Lambda_r(\ov D)$,  $f$ has an extension $\tilde f\in\Lambda_r(\rr^n)$ so that $|\tilde f|_{\Lambda_r}\leq 2|f|_{\Lambda_r}$. Take a decomposition $\tilde f=\sum f_j$ satisfying \re{nabi}. By \re{eSE}, we have  $\frak Ef=\sum \frak E(f_j|_{D})$ and the decomposition satisfies \re{nabi}, i.e. $| \mathfrak E(f_j|_{D})|_{L_k^\infty}\leq C_k|f_j|_{L_k^\infty(D)}\leq 2C_k |f|_{\Lambda_r}2^{j(r-k)}$. This shows that $|\mathfrak E  f|_{\Lambda_r(\rr^n)}\leq C'|f|_{\Lambda_r(\ov D)}$. We have verified \re{EL} for $D=D_i$.

We now verify \re{EH} for $D_i$.
 When $r$ is an integer and $f\in C^r(\ov D)$,  we need only to verify, by \re{eSE}, that $\pd^r\mathfrak Ef$ is continuous in $\rr^n$.
 And if $\all=r-[r]>0$, we need to estimate the $C^\all$ norm of $g=\pd_x^{[r]}\mathfrak Ef$.
Let us consider the case $[r]=0$ first. The continuity of $\mathfrak Ef$ follows  from the continuity of $f$, $|f|_{L^\infty}<\infty$,  the rapidly decreasing property of $\psi$, and
  \begin{equation*}
\mathfrak E f(x)-f(x',x_n+\delta^*(x))=\int_1^\infty R_0f(x,\la)\psi(\la)\, d\la,
  \end{equation*}
where $
 R_0(x,\la):=
 f(x',x_n+\la\delta^*(x))-f(x',x_n+\delta^*(x)).$  To estimate the H\"older ratio at two points $u,v$ in $\rr^n$, we may assume that $u,v$ are in $D^c$. Let $d=|u-v|$. Since $\del^*$ vanishes on $\pd D$ and $\pd_x\del^*$ is bounded in $\ov D^c$, then by connecting $u,v$ in the line segment we show that $|\del^*(u)-\del^*(v)|\leq C|u-v|$.
 Computing the H\"older ratio of each term in $R_0$, we obtain
 $$
 |R_0(u,\la)-R_0(v,\la)|\leq C_1\la^\all |f|_\all d^\all.
 $$
 We have verified \re{EH} for $0\leq r<1$.

For $k=[r]>0$,
 a $k$-th derivative $\pd^k\mathfrak Ef$ is a sum of $\mathfrak E\pd^kf$ and terms  of the form
\gan I(x)=\pd^{1+\ell_1}\del^*(x)\cdots\pd^{1+\ell_i}\del^*(x)\int_1^\infty\la^j\pd^{j}f(x',x_{n}+\la \del ^*(x))\psi(\la) \,d\la,
\end{gather*}
with $j+\ell_1+\dots+\ell_i\leq k$,  $ j>0$ and $i\leq j$. By the result for $[r]=0$, we know that $\mathfrak E\pd^kf$ is continuous and has the desired estimate in $|\cdot|_\all$ norm.
 With $j>0$ and \re{intp1}, we   subtract the Taylor polynomial of $\pd^{j}f(x',x_{n}+\la \del ^*(x))$ in $\la$
 of degree $k-j$ about $\la=1$ from $\pd^{j}f(x',x_{n}+\la \del ^*(x))$ and apply a Taylor remainder formula to express $I(x)$ as a linear combination of
$$
\tilde I(x)=\eta(x)\int_1^\infty\la^{j}R_kf(x,\la)\, \psi(\la) \,d\la,
$$
where $\eta(x):=\del^*(x)^{k-j}\pd^{1+\ell_1}\del^*(x)\cdots\pd^{1+\ell_i}\del^*(x)$ and
$$
R_kf(x,\la):= \int_1^\la(\la-\theta)^{k-j}\left\{\pd^{k}f(x',x_{n}+\theta\del ^*(x)) -\pd^{k}f(x',x_{n}+ \del ^*(x))\right\}\, d\theta.
$$
By \re{sDel}, we have
$
|\eta(x)|\leq C_0. 
$
It is now clear that $\tilde I$ and hence $I$ is continuous in $\ov D^c$ and vanishes in $\pd D$. Therefore $\mathfrak Ef\in C^{[r]}(\rr^n)$.

For the H\oo lder ratio, let $\all=r-[r]$. Take two points $x,x'\in\rr^n\setminus D$. If $x'\in\pd D$, we have $R_kf(x',\la)=0$ and
$$
|\tilde I(x)-\tilde I(x')|\leq C_0'|f|_{\ov D;r}\del^*(x)^{\all}\leq C_0'|f|_{\ov D;r}|x-x'|^{\all}.
$$
Let $d=|x-x'|$ and let $d_L$ be the distance from  $\pd D$ to the
 line segment $L$ connecting $x,x'$.
We consider two cases: (i)   $d_L\leq d$;
(ii) $d_L>d$.  In the first case, we take $x''\in\pd D$ with distance at most $d$ from $L$. Then $|x-x''|\leq2d$ and $|x'-x''|\leq2d$. We get
\aln
|\tilde I(x)-\tilde I(x')|&\leq |\tilde I(x)-\tilde I(x'')|+|\tilde I(x'')-\tilde I(x')|\\
&\leq C|f|_{a}(|x-x''|^{\all}+|x'-x''|^{\all})\leq C'|f|_ad^{\all}.\end{align*}
In the second case, we have
$
|R_kf(x,\la)-R_kf(x',\la)|\leq C\la^{k+1}|f|_{r}|x-x'|^{\all}.
$
Thus $\tilde \eta(x):=\int_1^\infty\la^{j}R_kf(x,\la)\, \psi(\la) \,d\la$ satisfies
\eq{ettet}
|\eta(x)(\tilde\eta(x)-\tilde\eta(x'))|\leq C|f|_r |x-x'|^{\all}.
\eeq
By \re{sDel}, we   have
$ |\pd\eta(x)|\leq C_1d(x)^{-1}.
$
Thus, $|\eta(x)-\eta(x')|\leq \sup_{\zeta\in L}|\pd_\zeta\eta||x-x'|\leq Cd_L^{-1}|x-x'|$ and $|\tilde\eta(x')|\leq C|f|_rd(x')^{\all}$, and we obtain
\eq{tetet}
|\tilde\eta(x')(\eta(x)-\eta(x'))|\leq C|f|_rd(x')^{\all}d_L^{-1}|x-x'|\leq C|f|_rd(x')^{\all}d_L^{-\all}|x-x'|^\all.
\eeq
Furthermore, $d_L\geq d(x')-|x-x'|\geq d(x')-d_L$. We simplify \re{tetet} and combine it with \re{ettet} to conclude   $|\tilde I(x)-\tilde I(x')|\leq C_r|f|_{r}|x-x'|^\all$, for the second case.

 Therefore, we have
verified \re{EH} for each $D_i$.
In the general case, we can verify that  the linear extension operator $\mathfrak E$ defined in~\ci{St70}*{p.~191, formula (31)}  satisfies \re{EL}-\re{EH}. We leave the details to the reader.
\end{proof}

\subsection{Real interpolation.}
We recall   some basic results of the real interpolation via the $K$-method of Peetre
from   Butzer-Berens~\ci{BB67}*{Sect. 3.2, p.~165}.  Let $X_0,X_1$ be two Banach spaces embedded continuously in a linear Hausdorff space $\cL X$. Define
$$|f|_{X_0\cap X_1}=\max\{|f|_{X_0},|f|_{X_1}\}, \quad |f|_{X_0+X_1}=\inf_{f=f_0+f_1}(|f_0|_{X_0}+|f_1|_{X_1}).$$
The $(X_0,X_1)$ is called an interpolation pair of Banach spaces in $\cL X$.  Define
$$
K(t,f;X_0,X_1)=\inf_{f=f_0+f_1}(|f_0|_{X_0}+t|f_1|_{X_1}),\quad t>0, \quad f\in X_0+X_1.
$$
Let $\theta\in(0,1)$. By $f\in X_{\theta,\infty;K}$, we mean that
$$
|f|_{\theta;X_0,X_1}:=\sup_{t>0} \{t^{-\theta}K(t,f;X_0,X_1)\}<\infty.
$$
Then $(X_0,X_1)_{\theta}:=X_{\theta,\infty;K}$ with norm $|\cdot|_{\theta;X_0,X_1}$ is a Banach space, while $X_0\cap X_1\subset (X_0,X_1)_{\theta} \subset X_0+X_1$ are continuous embeddings.
Following Triebel~\ci{Tr95}*{Sect.\!  2.7, p.p.~200-202}, we let $\mathcal{C}^r(\rr^n)$ be the closure of the space of rapidly decreasing functions in $\rr^n$ in $\Lambda_r(\rr^n)$.  Then by \cite{Tr95}*{Thm.\! 1, p.\! 201; Thm.\! (g), p.\! 50} we have
\eq{interL}
(\mathcal C^{r_0}(\rr^n),\mathcal C^{r_1}(\rr^n))_{\theta}=\mathcal C^{(1-\theta)r_0+\theta r_1}(\rr^n), \quad0<\theta<1, \
0<r_0<r_1<\infty
\eeq
in equivalent norms.
Let  $(Y_0,Y_1)$ be an interpolation couple of Banach spaces continuously embedded in a linear Hausdorff space   $\cL Y$. If $T\colon \cL X\to\cL Y$ is   linear, and if
$$
\|Tf_i\|_{Y_i}\leq M_i\|f_i\|_{X_i}, \quad i=0,1
$$
then $\|Tf\|_{(Y_0,Y_1)_{\theta}}\leq M_0^{1-\theta}M_1^\theta\|f\|_{(X_0,X_1)_{\theta}}$; see~\ci{BB67}*{Thm.\! 3.2.23, p.\! 180} or~\ci{Tr95}*{p.\! 26}.

In summary, we can apply  the following:
\pr{interpolation}Let $\mathcal C^r=\mathcal C^{r}(\rr^n)$.
Let $a_i,b_i$ be positive real numbers satisfying $a_0<a_1$ and $b_0<b_1$.  Let $T\colon\mathcal C^{a_0}\to\mathcal C^{b_0}$ be a linear operator such that
$
|Tf|_{\mathcal C^{b_i}}\leq M_i|f|_{\mathcal C^{a_i}},$ for $i=0,1.
$
Then in equivalent norms,
$|Tf|_{\mathcal C^{b_\theta}}\leq C_{r,b,\theta}M_0^{1-\theta}M_1^{\theta}|f|_{\mathcal C^{a_\theta}}$
 for $0<\theta<1$, $a_\theta=(1-\theta)a_0+\theta a_1$, and $b_\theta=(1-\theta)b_0+\theta b_1$.
\epr

\subsection{$\db$ solutions for the top type.}
As an application of the extension and interpolation, we   estimate a $\db$ solution for forms of type $(0,n)$. Let $C^r_{(0,q)}(\ov D)$ (resp. $\Lambda^r_{(0,q)}(\ov D)$) be the set of $(0,q)$-forms in $D$ of which the coefficients are in $C^r(\ov D)$ (resp. $\Lambda_r(\ov D)$). It seems that the following statement has not appeared in the  literature.
\pr{qisn}Let $D$ be a bounded domain in $\cc^n$.
\bppp
\item  Suppose that any two points $p,q$ in $D$ can be joined by a broken line segment $\gaa$ in $\ov D$ of length at most $ L|p-q|$, while $\gaa\cap\pd D$ is a finite set.
 For each $r\in(0,\infty)\setminus\nn$, there is a  linear map $T_r\colon C_{0,n}^r(\ov D)\to C_{0,n-1}^{r+1}(\cc^n)$, which depends on $r$,  so that $\db T_r\var=\var$ in $D$ and $|T_r\var|_{\cc^n;r+1}\leq C_r(D)|\var|_{\ov D;r}$.
\item Assume that
 $\pd D$ is Lipschitz. There is a linear operator $S\colon C_{0,n}(\ov D)\to C_{0,n-1}(\cc^n)$ so that $\db S\var=\var$ and $|S\var|_{\Lambda_{(0,n-1)}^{r+1}(\cc^n)}\leq C_r(D)|\var|_{\Lambda_{(0,n)}^{r}(\ov D)}$ with $C_r(D)<\infty$ for all $r\in(0,\infty)$.
\eppp
\epr
\begin{proof} $(i)$ We apply the Whitney extension $E_r $ for $\ov D$ via \rl{estR} and \rp{Whit}. Fix an open ball $B$ containing $\ov D$. By the Leray-Koppleman solution operator $T_n$ for $B$ and estimate in~\ci{We89}, we get the conclusion.

$(ii)$ Let $E\colon C^0(\ov D)\to C^0(\cc^n)$ be the bounded linear Stein extension. Thus
$E\colon\Lambda_a(\ov D)\to\Lambda_a(\cc^n)$ is bounded for all $a\in(0,\infty)$. We first consider the case of a non-integer $r$.
We have $\var=fd\ov z_1\wedge\cdots\wedge d\ov z_n$. By \re{EL}, we may assume that $D$ is relatively compact in a ball $B_0$. Replacing $\var$ by $E\var$, we may assume that   $\var\in \Lambda_{r}(\cc^n)$. Take a sequence $\var_j\in C^\infty(\cc^n)$ satisfying $\var_j\to \var$ in $L^\infty(D)$. Then we have $T_n\var_j\to T_n\var$ in $L^\infty(D)$.
Since $\db T_n\var_j=\var_j$, we get $\db T_n\var=\var$ in the sense of distributions.   By \re{Wea+1}, we get
$
|T_n\var|_{\ov D;r+1}\leq C_r|\var|_{B_0;r}.
$
By \re{EL} again, we conclude that \eq{boundE}
|ET_n\var|_{\cc^n;r+1}\leq C_r|\var|_{B_0;r}, \quad r>0.
\eeq
When $r$ is a positive integer,  the estimate follows from interpolation by \rp{interpolation} as follows.  We consider a linear operator
$$
\tilde T_n:=\chi ET_n\colon \mathcal C_{(0,n)}(\cc^n)\to \cL C_{(0,n-1)}(\cc^n),
$$
where $\chi\in C^\infty_0(\cc^n)$ and $\chi=1$ in $B_r$. By \re{boundE}, we have $|\tilde T_n\var|_{\Lambda_{r+1}}\leq M_r|\var|_{\Lambda_{r}}$ for $r\in(0,\infty)\setminus\nn$. By \rp{interpolation}, we get the same estimate for all positive integer $r$.
\end{proof}

\begin{rem}\label{rem1}The constant  $C_r(D)$ in \rp{qisn} (i) depends on $L$ and the diameter of $D$ and the $C_r(D)$ in (ii) depends on the constants   $\e,M, N$ in  Definition~\ref{mins}, as well as the diameter of $D$ by \rp{estL}.
\end{rem}

\section{Estimates for the homotopy operators}
\label{sect:4}
 \setcounter{thm}{0}\setcounter{equation}{0}

In this section we first introduce  a regularized Leray map to study strictly pseudoconvex domains with low regularity. The main estimates are derived under the assumption of the existence of a {\it regularized}  Henkin-Ram{\'{\i}}rez  function for the homotopy operators $H_q$.

\begin{defn}Let $D$ be a bounded domain of class $C^2$ and define
$$
D_\del=\{z\in\cc^n\colon\dist(z,\ov D)<\del\},\quad D_{-\del}=\{z\in D\colon\dist(z,\pd D)>\del\}, \quad \del>0.
$$
 We say that
$W$
is a {\it regularized} Leray mapping in $D_\del\times(D_\del\setminus D_{-\del})$, if  for some positive number $\del$ the following hold
\bppp
\item
$
W\colon D_{\del}\times( D_\del\setminus D_{-\del})\to\cc^n
$ is a $C^1$ mapping,  and  $W(z,\zeta)$ is holomorphic in $z\in D_{\del}$.
\item $W(z,\zeta)\cdot(\zeta-z)\neq0$ for $z\in D$ and $\zeta\in
D_{\del}\setminus D$.
\item For each $z\in D_{\del}$, we have $W(z,\cdot)\in C^1(\ov D\setminus D_{-\del})$ and
$$
|\pd_\zeta^iW(z,\zeta)|\leq C_i|W(z,\cdot)|_{\ov D,1}(1+\dist^{1-i}(\zeta,D)), \quad \zeta\in D_{\del}\setminus\ov D, \quad 0\leq i<\infty.
$$
\eppp
\end{defn}
The first two properties are the standard requirements for the Leray maps. The third   is new.
The existence of a regularized $C^2$ defining function for a domain with $C^2$ boundary is proved in \rl{regdef}. The Whitney extension of a strictly convex function $\rho$ in $\ov D$ remains strictly convex in a neighborhood of $\ov D$. Therefore, we have the following.
\begin{exmp}
Let $D$ be defined by $\rho_0<0$  in ${\cL U} $ with $\ov D\subset {\cL U} $. Suppose that $\rho_0$ is a $C^2$ {\it strictly} convex function in ${\cL U} $. Let $\rho$ be a Whitney extension of $\rho_0|_{\ov D}$ as in \rla{regdef}. Then
$W( z,\zeta )=(\rho_{\zeta_1},\dots, \rho_{\zeta_n})$ is a regularized   Leray mapping.
\end{exmp}

%
%

We now derive our main estimates. Recall the homotopy operator
$$
H_q\var=\int_{\cL U}\Om_{0,q-1}^0\wedge E\var+\int_{\cL U\setminus D}\Om_{0,q-1}^{01}\wedge[\db,E]\var.
$$
Here $\cL U=D_\del$. The first term is estimated by \rp{estL}, gaining one derivative in a H\"older space. We now estimate the second term for $z\in D$:
 \eq{dbE}
 \int_{\cL U\setminus D}\Om_{0,q}^{01}(z,\zeta)\wedge[\db,E]\var(\zeta).
 \eeq
 From now on, we   take $g^0(z,\zeta)=\ov z-\ov \zeta$ and $g^1(z,\zeta)=W(z,\zeta)$. We require that $W$ is a regularized Leray map.

We will denote by $\pd_z^k$ a derivative of order $k$ in $(z,\ov z)$, and by $N_k(\zeta-z)$
  a monomial in $\zeta-z, \ov\zeta-\ov z$ of degree $k$. Let $A(w)$ denote a polynomial in $w,\ov w$, where $N_k$ and $ A$ may differ when they recur.
  We can write
 \re{dbE} as a linear combination of
 \ga\label{defnKf}
Kf(z):=\int_{\cL U\setminus D} f(\zeta )\f{A( \hat\pd_\zeta W(z,\zeta), z,\zeta )N_{1}(\zeta-z )}
 {\Phi^{n-\ell}(z,\zeta)|\zeta -z |^{2\ell}}\, dV(\zeta), \quad 1\leq\ell<n,\\
 \Phi(z,\zeta)=W(z,\zeta)\cdot(\zeta-z),
 \label{defnKf+}
 \end{gather}
 where $f$ is a coefficient of the form $[\db,E]\var$. In particular $f$ vanishes on $\ov D$. Here $\hat\pd_\zeta W$ denotes $W$ and its first-order $\zeta$ derivatives.

To derive our main estimates, we start with the following   lemma.
\le{a+.5}Let $\beta\geq 0$,  $\all\geq0$,  and let $0<\del<1/2$.
\bppp
\item If   $\all<\beta+1/2$, then
 $
\int_{0}^1\int_{0}^1\f{s^{\all+1}\, dt \, ds}{(\del+s+t^2)^{3+\beta} }\leq C\del^{\all-\f{1}{2}-\beta}.
$
\item
$\int^{s=2\del}_{s=\del}\int_{0}^1\f{s^{\all+1}\, dt \, ds}{(s+t^2)^{1+\beta} }\leq C\del^{\all-\beta+3/2}.$
\eppp
\ele
\begin{proof}  (i). We divide  $[0,1]\times[0,1]$ in the $(s,t)$-plane  in three regions
$$
P\colon \del+s\geq t, \quad Q\colon \del+s\leq t^2,\quad R\colon t\geq \del+s\geq t^2.
$$
  The integral in $P$ is bounded above by
\gan
 \int_{s=0}^1\int_{t=0}^{\del+s}\f{s^{\all+1} \, dt\, ds}{(\del+s)^{3+\beta}}
\leq \int_{0}^1\f{s^{\all+1} \, ds}{(\del+s)^{2+\beta}}\leq \int_{0}^1(\del+s)^{\all-1-\beta}\, ds,\end{gather*}
which is less than $C\del^{\all-\beta-1/8}$. In $Q$, it is bounded  by
\gan
 \int_{s=0}^1\int_{t=\sqrt{\del+s}}^1\f{s^{\all+1} \, dt\, ds}{t^{6+2\beta}}\leq   \int_{s=0}^1\f{s^{\all+1}}{(\del+s)^{\beta+5/2}} \, ds,\end{gather*}
 which is less than $C\del^{\all-\beta-1/2}$.
 In $R$, it has a similar  bound as
\gan
 \int_{s=0}^1\int_{t=0}^{\sqrt{\del+s}} \f{s^{\all+1} \, dt\, ds}{(\del+s)^{3+ \beta}}\leq   \int_{s=0}^1\f{s^{\all+1} \, ds}{(\del+s)^{\beta+5/2}}.
 \end{gather*}

 (ii). We divide the domain $[\del,2\del]\times[0,1]$  in three regions
 $$
P\colon  s\geq t; \quad Q\colon  s\leq t^2;\quad R\colon t\geq s\geq t^2, \quad t\leq \del^{1/2}.
$$
The integral in $P$ is bounded above by
$
\int_\del^{2\del}s^{\all-\beta}\int_0^sdt\, ds\leq C\del^{\all-\beta+2}.
$
In $Q$, it is bounded   by
$
\int_\del^{2\del}s^{\all+1}\int^{1}_{\sqrt s}t^{-2(\beta+1)}\,   dt\, ds\leq C\del^{\all-\beta+3/2}.
$
In $R$, it is bounded   by
$
\int_\del^{2\del}s^{\all-\beta}\int_{0}^{\sqrt s}\, dt\, ds\leq C\del^{\all-\beta+3/2}.
$
 \qedhere
\end{proof}

\pr{estK} 
Let $1\leq r<\infty$ and $\all=r-[r]$. Let $D$ be a strictly pseudoconvex domain defined by $\rho<0$ with $\rho\in C^2(\cL U)$ and $\ov D\subset\cL U$. Suppose that $\pd \rho\neq0$ in $\pd D$.
Let $W$ be a regularized Leray mapping of $D$
in $D_\del\times(D_\del\setminus D_{-\del})$ and let $\Phi, Kf$ be defined by \rea{defnKf}-\rea{defnKf+}.
 Assume that there is a finite open covering  $\{\om_1,\dots, \om_m\}$ of $\pd D$ and $C^1$ coordinate maps $\Psi_i\colon \zeta\to (s,t)=(s_1,s_2,t_3,\dots, t_{2n})$ defined in  $\om_i$  such that $s_1=\rho(\zeta)$ and
 for $ z\in\om_i\cap D,\zeta\in\om_i\setminus D$
  \ga\label{LbPhi}
  |\Phi(z,\zeta)|\geq  c_*(d(z)+ s_1(\zeta)+|s_2(z,\zeta)|+|t(z,\zeta)|^2),\\
 |\Phi(z,\zeta)|\geq c_*|\zeta-z|^2,\quad |\zeta-z|\geq c_* |(s_2,t)(z,\zeta)|,
 \label{LbPhi+}
  \end{gather}
  where $c_*>0$ is a constant.
  Suppose that 
  $f$ vanishes in $D$.
  Then the following hold:
  \bppp
  \item
  Suppose that $f\in C_0^{r-1}(\cL U)$.
   If $r\geq1$, then for $z\in D$ and $d(z)=\dist(z,\pd D)$
\ga  \label{aneqhalf}
|\pd_z^{[r]+1}Kf  | \leq C_r  d(z)^{\all-1/2}\|f\|_{\cL U;r-1}, \quad 0\leq\all<1/2,\\
|\pd_z^{[r]+2}Kf  | \leq C_r   d(z)^{\all-3/2}\|f\|_{\cL U;r-1}, \quad 1/2\leq\all<1.
\label{aeqhalf+}
\end{gather}
In particular,    $\|Kf  \|_{\ov D,r+1/2}\leq C_r\|f\|_{  \cL U,r-1}$ for $\all\neq1/2$.
    \item
 Suppose that $f\in \Lambda_0^{r-1}(\cL U)$ and  $r>1$. Then
 $\|Kf  \|_{\Lambda_{r+1/2}(\ov D)}\leq C_r\|f\|_{\Lambda_{r-1}(\cL U)}$.
\eppp
Here $C_r$ depends on $r$, $c_*$,  $\sup_{z\in D_\del}|W(z,\cdot)|_{D_\del\setminus D_{-\del};1}$, and $C^1$ norms of $\Psi_i$ and the sup norms of $(\det\pd_\zeta\Psi_i)^{-1}$.
\epr
\begin{proof} By the assumption, we have $\pd D\in C^1$.
We have  $\Phi(z,\zeta)\neq0$,  for $z\in D$ and $\zeta\in {\cL U} \setminus D$.
 The latter contains
the support
of $f$.  We will first consider the case  $f\in C^{r-1}$. We will  distribute  the first $[r]-1$ derivatives of $Kf(z)$ directly to the integrand when $[r]>1$. We will then apply   the integration by parts in $\zeta$ variables to derive a new formula.

$(i)$ By the assumption, we have $f\in C_0^{r-1}(\cL U)$ and $f=0$ in $D$. We may assume that $\cL U=D_\del$.
We   write
$
\pd_z^{[r]-1}\{Kf(z )\}$ as a linear combination of  $  K_1f(z)$ with
\ga\label{K1gx}
K_1f(z):=\int_{\cL U\setminus \ov D} f_1(z,\zeta )\f{
N_{1-\mu_0+\mu_2}(\zeta -z )}
 {\Phi^{n-j+\mu_1}(z,\zeta)|\zeta -z |^{2j+2\mu_2}}\, dV(\zeta),\\
 f_1(z,\zeta)=f(\zeta)A_1 ( W_1(z,\zeta), z,\zeta ),\quad
 W_1(z,\zeta)=(\hat \pd_\zeta W(z,\zeta),\pd_z^{k_0}\hat \pd_\zeta W(z,\zeta)),\\
 \label{aaa}
\mu_0+\mu_1+\mu_2\leq [r]-1, \quad 1-\mu_0+\mu_2\geq0,\quad k_0\leq[r]-1,
 \end{gather}
 where $A_1$ is a polynomial.

Let us
 explain how $\pd D\in C^2$ suffices the estimation.  Since $\hat\pd_\zeta W(z,\zeta)$ is holomorphic in $z$,  its $z$-derivatives in a suitable neighborhood of $D_{\del/2}$
  can be estimated by the sup-norm of $\hat\pd_\zeta W$ in $D_\del$ by using the Cauchy formula.
In $\cL U\setminus\ov D$, the integrand in $K_1f$ is smooth in $z$. As $\zeta\in\cL U\setminus\ov D$ approaches $\pd D$, the rate of growth of a $\zeta$-derivative of $W_1$ is bounded by a precise  negative power of $ d(\zeta)$. The latter can be dominated by the order of vanishing of $f(\zeta)$ along $\pd D$. Let us record the estimate
\ga\label{pdzW}
|\pd_\zeta^i\pd_z^jW_1(z,\zeta)|\leq C_{i+j}(|W|_{1})
d(\zeta)^{-i},  \\
C_{k}(|W|_1):=C_k(\sup_{z\in D_{\del}}|W(z,\cdot)|_{\ov D_\del,1}),
\end{gather}
for $  \zeta\in \cL U\setminus\ov D, z \in D_{\del/2}$ and $i,j=0,1,\dots$.
 Our argument
relies essentially on  that   $f(\zeta)$ is independent of $z$.


 We now provide the details of the proof. We will use  integration by parts as in Elgueta~\ci{El80},  Ahern-Schneider~\ci{AS79}, and Lieb-Range~\ci{LR80}
   to  reduce the exponent
    of $\Phi$  to the original
    $n-j. $
In our case,  the integration by parts needs to be carried out by using \rl{regstokes},  since $W_1(z,\zeta)$  is $C^\infty$ in $\zeta\in D_\del\setminus \ov D$ and it is, however,  merely $C^1$ in $\zeta\in D^c$.
     To this end, we write \re{K1gx}   as
 \ga
  K_1  f(z):=\int_{\cL U\setminus \ov D}\f{   h_1( z,\zeta ) }{\Phi^{n-j+\mu_1}(z,\zeta)}\, dV(\zeta),
\nonumber 
\\
\intertext{with}
 h_1( z,\zeta )=  f_1( z,\zeta )\f{  N_{1-\mu_0+\mu_2}(\zeta-z)}
 { |\zeta-z|^{2j+2\mu_2}}.
\nonumber 
 \end{gather}
 Using a
 partition of  unity in $\zeta$ space and replacing $f$ by $\chi f$ for a $C^\infty$ cut-off function, we may assume that
 \ga
\nonumber 
 \supp f\subset B_0\setminus D,\quad
 u( z,\zeta ):=\pd_{\zeta_{i^*}}\Phi(z,\zeta)\neq0
\nonumber 
\end{gather}
  for some $i^*$. Here $B_0$ is a small open set in $\cL U$ containing a given $\zeta_0\in\pd D$.  Recall that $\pd D\in C^1$. We have for $\zeta\in B_0\setminus \ov D$
\gan 
|\pd_\zeta^iW(z,\zeta)|\leq C_i(|W|_1) (1+d(\zeta)^{1-i}),\quad i=0,1,\dots,\\
|\pd_\zeta^{i}W_1(z,\zeta)|+\left|\pd_\zeta^i\f{1}{u(z,\zeta)}\right|\leq C_i(|W|_1) d(\zeta)^{-i}, \quad i=0,1,\dots.
\end{gather*}
Up to a constant multiple, we rewrite $K_1f$ as
$$
K_1f=\int_{\cc^n\setminus \ov D} u(z,\zeta)^{-1}  h_1( z,\zeta ) \pd_{\zeta_{i^*}}\Phi^{-(n-j+\mu_1-1)}(z,\zeta)\, dV(\zeta).
$$
Since $f\in C^{r-1}$ vanishes identically in $\ov D$, then $\pd^if$ vanishes in $\pd D$ for $i\leq [r]-1$.  Thus, by Taylor's theorem,
\eq{estdf}
|\pd^if(\zeta)|\leq C_i|f|_{D_\del;r-1}d(\zeta)^{r-1-i}, \quad \zeta\in D_\del, \quad 0\leq i\leq[r]-1.
\eeq
Suppose that $[r]>1$. Fix $z\in D$. Thus $|z-\zeta|^{-2j-2\mu_2}$ is $C^\infty$ in $\zeta\in \cc^n\setminus D$. Recall that $f_1(z,\zeta)=f(\zeta)A_1 ( W_1(z,\zeta), z,\zeta )$. Using \re{pdzW} and \re{estdf},  a straightforward computation shows that
\ga\label{u-1h}
|\pd_\zeta^{i}((u^{-1}h_1)(z,\zeta))|\leq C_{i,j}(z)C_{i+j}(|W|_1)|f|_{D_\del;r-1}
d(\zeta)^{r-1-i},
\end{gather}
for $\zeta\in D_\del\setminus\ov D$ and $j\in\nn$. Here $C_{i,j}(z)<\infty$ because $z\in D$.    In particular, this allows us to apply the integration by parts to transform $K_1f$.

When $[r]>1$, we apply Stokes' theorem via \rl{regstokes} in which $S(\zeta)=u(z,\zeta)^{-1}h_1(z,\zeta)$,  $B(\zeta)=\Phi(z,\zeta)^{-(n-j+\mu_1-1)}$, $m=[r]-1$,  and $b=0$.  (Recall that we fix $z\in D$. Thus $B$ is  $C^1$ in $D^c$ and $C^\infty$ in its interior.)
Up to a constant multiple, we have
  $$
    K_1 f(z)=\int_{\cc^n\setminus \ov D} \f{ \pd_{\zeta_{i^*}} \{u( z,\zeta )^{-1}  h_1( z,\zeta )\}}{\Phi^{n-j+\mu_1-1}(z,\zeta)}\, dV(\zeta),
  \quad\forall z\in  D. $$
 We also have, up to a constant multiple depending on $z,|W|_1$,
\aln
&|\pd_\zeta^i \left\{ u(z,\zeta)^{-1} (\pd_{\zeta_{i^*}}\circ u( z,\zeta )^{-1})^{\ell}\{ h_1( z,\zeta )\}\right\}|  \leq  
|f|_{D_\del;r-1} d(\zeta)^{r-1-\ell-i},\\
&|\pd_\zeta^i\Phi^{-(n-j+\mu_1)+\ell}|\leq C_{i,\ell}(z) C_{i+\ell}(|W|_1)(1+d(\zeta)^{1-i}).
\end{align*}
If $\mu_1-\ell>0$, we have $[r]-1-\ell>[r]-1-\mu_1\geq0$ by \re{aaa}.  Applying  integration by parts $\mu_1$ times via \rl{regstokes} till  $\mu_1-\ell=0$, we obtain
 $$
  K_1  f(x):=\int_{\cc^n\setminus \ov  D} \f{ h_2( z,\zeta )}{\Phi^{n-j}(z,\zeta)}\, dV(\zeta),
  \quad\forall z\in   D$$
 with
 $$ 
  h_2( z,\zeta ):=
 (\pd_{\zeta_{i^*}}\circ u( z,\zeta )^{-1})^{\mu_1}\{   h_1( z,\zeta )\}.
$$ 
By the product and quotient rules, we can write $ h_2( z,\zeta )$ as a linear combination of
\ga\label{g4xix}
 h_3( z,\zeta )=f_2( z,\zeta )\widetilde N_\la(\zeta-z),
\end{gather}
with
\ga
\widetilde N_\la(\zeta-z)=\f{  N_{1-\mu_0+\mu_2-\nu_1+\nu_2}( \zeta-z)}
 { |\zeta-z|^{2j+2\mu_2+2\nu_2}},\\
f_2( z,\zeta )=A_2\bigl(\widetilde W_1(z,\zeta), z,\zeta \bigr)\,\pd_\zeta^{\nu_i}\widetilde W_1(z,\zeta)\cdots \pd_\zeta^{\nu_4}\widetilde  W_1(z,\zeta)\,
\pd^{\nu_{3}}_\zeta f.
\label{f5}
\end{gather}
Furthermore, each $\widetilde  W_1$ is  one of $\hat \pd_\zeta W(z,\zeta),\pd_z^{k_0}\hat \pd_\zeta W(z,\zeta),u_1^{-1}(z,\zeta)$. And
\ga
\label{g5xi}
 \nu_1+\dots+\nu_i\leq\mu_1, \quad 1-\mu_0+\mu_2-\nu_1+\nu_2\geq0,\\
 \label{Lamu} \la:=(1-\mu_0+\mu_2-\nu_1+\nu_2)-(2j+2\mu_2+2\nu_2).
\end{gather}
We have proved that $\pd_z^{k-1}Kf$ is a linear combination of
$K_1f$, while $K_1f$  is a linear combination of
\eq{aneqhalf++}
K_2f(z)=\int_{\cc^n\setminus\ov  D}f_2( z,\zeta )\f{\tilde N_\la (\zeta-z)}{\Phi^{n-j}( z,\zeta )}\, dV(\zeta).
\eeq
Since $f(\zeta)=0$ in $\ov D$, it is easy to see that $  K_2  f\in \cL C^\infty(D)$.

We want to estimate $K_2f(z)$ in terms of distance $ d(z)$.  To achieve the estimate that has the form \re{aneqhalf}, we need to count the numbers of derivatives in the expression of $f_2$. In \re{f5}, we have applied $\nu_4+\dots+\nu_i$ extra derivatives on $\widetilde W_1$. Set
$$
\la'=\nu_4+\cdots+\nu_i.
$$
Since $|\zeta-z|\geq d(\zeta)$ and $[r]-1- \nu_3\geq0$ by \re{aaa} and \re{g5xi},  we obtain
for  $\zeta \in\ov D^c$,
$$
|\pd_\zeta^{\nu_3}f|\leq C|f|_{r-1} d(\zeta)^{r-1- \nu_3}, \quad
|\pd_\zeta^{\nu_i}\widetilde W_1(z,\zeta)\cdots \pd_\zeta^{\nu_4}\widetilde  W_1(z,\zeta)| \leq   C(|W|_1) d(\zeta)^{ -\la'},
$$
where the last inequality follows from \re{pdzW}. Hence their $z$-derivatives can be estimated by $|W|_{D_\del,1}$.
Then we have proved that for $
z\in D$ and $\zeta\not\in \ov D$
\al\label{g5xiL}
|f_2( z,\zeta )|&\leq C(|W|_1) |f|_{r-1} d(\zeta)^{r-1-\la'-\nu_3}\\
 &\leq C(|W|_1)|f|_{r-1} d(\zeta)^{\all}|\zeta-z|^{[r]-1-\nu_3-\la'}\nonumber
\end{align}
by using
$
[r]-1-\nu_3-\la'\geq0$ and $ d(\zeta)\leq|\zeta-z|.
$ By the definition of $\tilde N_\la$, we have
\gan
|\widetilde N_\la(\zeta-z)|\leq|\zeta-z|^\la, \quad |h_2(z,\zeta)|\leq C(|W|_1)|f|_{r-1} d(\zeta)^\all|\zeta-z|^{[r]-1-\nu_3-\la'+\la}.
\end{gather*}
We have just estimated $h_2$. Since $f(\zeta)$ does not depend on $z$,
 the $z$-derivatives of $f_2(z,\zeta)$, given by \re{f5},  satisfy
\ga
\nonumber 
|\pd_z^\ell f_2( z,\zeta )|\leq C_i |f|_{r-1} d(\zeta)^{\all} |\zeta-z|^{[r]-1-\nu_3-\la'},
\quad  [r]-1-\nu_3-\la'\geq0,\\
\nonumber 
|\pd_z^\ell \widetilde N_\la(\zeta-z)|\leq C_\ell  |\zeta-z|^{\la-\ell },
\\
 |\pd_z^\ell\Phi^{-(n-j)}(z,\zeta)|\leq C_\ell (|W|_1) |\Phi^{-(n-j)-\ell}(z,\zeta)|.\nonumber
\end{gather}
We estimate $\tilde N_\la$ first. We have \aln
\la &= (1-\mu_0+\mu_2-\nu_1+\nu_2)-(2j+2\mu_2+2\nu_2)\\
\nonumber
&=1-2j -\mu_0-\mu_2-\nu_1-\nu_2\geq 1-2j -\mu_0-\mu_2- \mu_1+\nu_3+\la'\\
&\geq 2-2j-[r]+\nu_3+\la'
\nonumber
\end{align*}
by the first inequalities in \re{g5xi} and \re{aaa}. Thus $[r]-1-\nu_3+\la\geq1-2j$.
For $h_3$ given by \re{g4xix}, we have
\ga
 \label{estg4}
|\pd_z^\ell h_3( z,\zeta )|\leq C(|W|_1)|f|_{r-1} d(\zeta)^{\all}|\zeta-z|^{1-2j-\ell }.
\end{gather}

We have expressed $\pd^{[r]-1}Kf$ as a linear combination of $K_2f$
by exhausting  all derivatives of $f$. Let $z\in D$.   We want to show that
 \ga
 \label{a121}
 |\pd_{z}^2K_2  f(z)|\leq C(|W|_1)|f|_{r-1} d(z)^{-1+(\all+ 1/2)},\quad \all+1/2<1,\\
\label{a122}  |\pd_{z}^3K_2  f(z)|\leq C(|W|_1)|f|_{r-1} d(z)^{-1+(\all-1/2)},\quad  \all+1/2\geq1.
 \end{gather}

   For $\ell =2,3$, we compute $\pd_z^\ell K_2f$ by differentiating the integrand directly. The $\pd_z^2K_2f$
is  a sum of three kinds of terms
\gan
J_0f(z)=\int_{\cc^n\setminus \ov D} \f{f_2(z,\zeta)}{\Phi^{n-j}(z,\zeta)}\, \pd_z^2 \left\{\tilde N_\la(\zeta-z)\right\}\, dV(\zeta),\\
 J_1f(z) =\int_{\cc^n\setminus\ov  D} \f{f_2(z,\zeta)}{\Phi^{n-j+1}(z,\zeta)}\,\pd_z\left\{\tilde N_\la(\zeta-z)\right\}\, dV(\zeta),\\
 J_2f(z) =\int_{\cc^n\setminus\ov  D} \f{f_2(z,\zeta)}{\Phi^{n-j+2}(z,\zeta)} \left\{\tilde N_\la(\zeta-z)\right\}\, dV(\zeta),
\end{gather*}
where $f_2$ still has the form \re{f5}  while $\nu_4,\dots, \nu_i$ are unchanged.
 Therefore we obtain
\ga
\nonumber 
|J_0f(z)|\leq C(|W|_1)|f|_{r-1}\int_{\cL U\setminus\ov  D}\f{ d(\zeta)^\all}
{|\Phi(z,\zeta)|^{n-j}|\zeta-z|^{1+2j}}\, dV(\zeta),\\
\nonumber 
|J_1f(z)|\leq C(|W|_1)|f|_{r-1}\int_{\cL U\setminus\ov  D}\f{ d(\zeta)^\all}
{|\Phi(z,\zeta)|^{n-j+1}|\zeta-z|^{2j}}\, dV(\zeta),\\
|J_2f(z)|\leq C(|W|_1)|f|_{r-1}\int_{\cL U\setminus\ov  D}\f{ d(\zeta)^\all}
{|\Phi(z,\zeta)|^{n-j+2}|\zeta-z|^{-1+2j}}\, dV(\zeta).
\nonumber %
\end{gather}
Recall that $1\leq j< n$.
For $z\in D$ and $\zeta\not\in D$, we have $C'|\zeta-z|\geq|\Phi(z,\zeta)|\geq C|\zeta-z|^2$. Thus it suffices to  estimate the last integral  for  $j=n-1$. Set
$$
\widehat J_2(z):=\int_{\cL U\setminus\ov  D}\f{ d(\zeta)^\all}
{|\Phi(z,\zeta)|^{3}|\zeta-z|^{2n-3}}\, dV(\zeta).
$$

Fix $\zeta_0\in\pd D$ and a small neighborhood $\om_0$ of $\zeta_0$.  Let $z\in\om_0\cap D$ and $\zeta\in\om_0\setminus D$.   Note that $r(\zeta)\approx \dist(\zeta,\pd D)= d(\zeta)$. We now use the assumption that
$$
|\Phi(z,\zeta)|\geq  c_*( d(z)+s_1(\zeta)+|s_2(z,\zeta)|+|t(z,\zeta)|^2), \quad |\zeta-z|\geq c_* |(s_2,t)(z,\zeta)|.
$$
  We also have
$$
 d(\zeta)/C\leq  r(\zeta)=s_1(\zeta) \leq |(s_1(\zeta),s_2(z,\zeta))|, \quad\zeta\in D_\del\setminus\ov D.
 $$
   Using polar coordinates for $(s_1,s_2)\in\rr^2$ and $(t_3,\dots, t_{2n})\in\rr^{2n-2}$, we  obtain for $z\in\om_0$
$$
 \widehat J_2(z)\leq C\int_{s=1}^1\int_{t=0}^1\f{s^{\all+1} \,
ds\, dt}{( d(z)+s+t^2)^3},
$$
which is less than $C d(z)^{\all-1}$ by \rl{a+.5} in which   $\beta=0$ and $0\leq\all<1/2$. We have verified \re{a121}.

Consider now the case $1/2<\all<1$. This requires us to estimate $\pd_z^3K_2f$, which
is   a sum of terms
$$
\widetilde J_if(z):= \int_{\cc^n\setminus\ov  D} \f{f_2(z,\zeta)}{\Phi^{n-j+i}(z,\zeta)}\, \pd_z^{3-i} \left\{\widetilde N_\la(\zeta-z)\right\}\, dV(\zeta)
$$
for $i=0,1,2,3$.  The worst term is $\widetilde J_3f(z)$ with $j=n-1$ and $i=3$. We have
 \aln
 |\widetilde J_3f(z)|
&\leq  C'(|W|_1)|f|_{r-1}\int_{s=1}^1\int_{t=0}^1\f{s^{\all+1} \,
ds\, dt}{( d(z)+s+t^2)^4},
 \end{align*}
which is less than $C|f|_{r-1} d(z)^{\all-3/2}$ by \rl{a+.5} with $\beta=1$ and $1/2<\all<3/2$.

%
%
%
%
%
%

\medskip

$(ii)$. We now consider the estimate in the $\Lambda_{r+1/2}$ space.

\medskip

 {\bf Case 1, $\all\neq0,1/2$.}
Recall that $Kf$ is $C^\infty$ in $D$ and its derivatives on a compact subset of $D$ can be estimated easily by the sup norm of $f$.
When $\all\neq0,1/2$,  by the Hardy-Littlewood lemma for H\"older spaces we get the estimate in $(ii)$ from $(i)$ immediately. Note that the same argument by the Hardy-Littlewood also gives us $|Kf|_{\ov D,k+1/2}\leq C_k|f|_{\cL U,k-1}$
when $k$ is a positive integer, which is however a weaker version of $(ii)$  for the $\Lambda_{k+1/2}$ estimate.

\medskip

 {\bf Case 2, $\all=1/2$.}
In this case \re{a121} says that
$
|\pd_z^{3}K_2f|\leq C|f|_{r-1}\dist(z)^{-1}.
$
 We remark that if we have $\pd D\in C^\infty$,
then by a version of Hardy-Littlewood lemma (see~\cite{MS94}), we could  conclude that $K_2f\in\Lambda_{1}$. Since $\pd D$ is only $C^2$, We need
another proof for the case $\all=1/2$, by using the estimates in $(i)$.

In fact we will provide an argument that actually works for $0<\all<1$.  Let us show that
$$
|K_2f|_{\Lambda_{r+1/2}}\leq C_r(|W|_1) |f|_{\Lambda_{r-1}}.
$$
We may assume that the $f$ vanishes when $|(s_1,s_2)|>1$ or $|t|>1$. We consider a  dyadic decomposition with
 $A_k^+:=\{(s_1,s_2)\colon 2^{-k-1}< |(s_1,s_2)|<2^{-k+1},s_1\geq0\}$ for $k=1,2,\dots$. Take a partition of unity $\{\chi_k\}$ such that $\supp\chi_k$ is contained in $A_K^+$, $\sum_k\chi_k=1$ in $\cup_k A_k^+$, and
 $|\pd ^j\chi_k|\leq C_j2^{jk}$ for $j=0,1,\dots$.
 Set $K_2f=\sum_{k\geq1}g_k$  with
$$
g_k(z)=\int_{\cL U\setminus \ov D}\chi_k(s_1(\zeta),s_2(z,\zeta))f_2(z,\zeta)\f{\tilde N_\la (\zeta-z)}{\Phi^{n-j}( z,\zeta )}\, dV(\zeta), \quad z\in D.
$$
 By \re{g4xix} and \re{estg4}
we still have $|h_3( z,\zeta )|\leq |f|_{r-1} s_1^{\all}(\zeta)|\zeta-z|^{1-2j}$.
  Now $\pd_z^ig_k$ is a linear combination of
$$
I_{i,\ell}(z):=\int_{\cL U\setminus \ov D}\pd_z^{i-\ell}\Bigl\{\chi_k(s_1(\zeta),s_2(z,\zeta))\Bigr\}\pd_z^{\ell}\left\{f_2(z,\zeta)\f{\tilde N_\la (\zeta-z)}{\Phi^{n-j}( z,\zeta )}\right\}\, dV(\zeta), \quad z\in D
$$
for $\ell=0,1,\dots, i$.  Again the worst term occurs to $j=n-1$ and $\ell=i$. Thus,
\aln
&|I_{i,\ell}(z)|\leq \int_{A^+_k} \int_{t\in\rr^{2n-2},|t|<1}\f{C'(|W|_1)|f|_{r-1}\cdot 2^{(i-\ell)k} s_1^{\all}}{(s_1+|s_2|+|t|^2)^{1+\ell}(s_1+|s_2|+|t|)^{2n-3}}\,ds_1\,ds_2\, dt\\
&\quad\leq C(|W|_1)|f|_{r-1}
\int_{s=2^{-k-1}}^{ 2^{-k+1}} \int_{t=0}^1\f{s^{\all+1}2^{(i-\ell)k}}{(s+t^2)^{1+\ell}}\,ds\, dt\leq C'(|W|_1)|f|_{r-1} 2^{-k(\all-i+1/2)},
\end{align*}
where the last integral is estimated in two regions $s\leq t^2$ and $s\geq t^2$.
Now assertion $(ii)$ for $0<\all<1$ follows from Proposition~\ref{steinext} and Lemma~\ref{gs141}.

\medskip
 {\bf Case 3,  $r>1$  an integer.}
We will achieve the $\Lambda_{r+1/2}$  estimate by the real interpolation theory.
Fix $d\ov z^I$ with $|I|=q>0$ and fix $d\ov z^J$ with $|J|=q-1$. Let $\{\psi\}_J$ denote the coefficients of $d\ov z^J$ for a $(0,q-1)$-form $\psi$. Consider the linear mapping
$$
L_J\colon f\mapsto \left\{ \int_{\cL U\setminus\ov  D}\Om_{0,q}^{01}\wedge[\db,E](f\, d\ov z^I)\right\}_J.
$$
Assume that $r \geq2$ be an integer.
Let $E$ be the linear extension operator  for functions defined in $\ov D$, given  in \rp{steinext}.
For the interpolation theory to be applicable, it is crucial that there is no other restriction to $f$.
Using \re{a122},  we have
\aln
|EL_Jf|_{\cc^n;r-\e+\yt}\leq C_1 |L_Jf|_{\cc^n;r-\e+\yt}\leq C_1C_k(|W|_{1})|f|_{\cc^n;r-1-\e}.
\end{align*}
 Using \re{a121},  we have
 $
|EL_Jf|_{\cc^n;r+\e+\yt}\leq   C_1C_{k+1}(|W|_{1})|f|_{\cc^n;r-1+\e}.
$
The estimate follows from  interpolation via  \rp{interpolation}.
 The assertion $(ii)$ is proved.
 \end{proof}

\begin{rem} In connection with Question $1$ in the introduction,
one can approximate   $\var\in\Lambda_1(\cc^n)$ by bounded $C^1$ forms $\var_j$ in $\cc^n$, which converges in the sup norm to $\var$. However, we  do not have  a useful limit $u$ of $H_q\var_j$ as $j\to\infty$, in order
to conclude that $\db u=\var$.
\end{rem}

We now turn to the estimate of holomorphic projection $H_0$. The analogous estimate for the boundary operator in \re{lerayf-} is in Ahern-Schneider~\ci{AS79}, where $\pd D\in C^\infty$ is used. We need to restrict  to $r>1$, requiring $\pd D\in C^2$ only.

\le{a+.5+}Let $0\leq\all<1$,  $0<\delta <1/2$, and $n\geq2$. Then
 $$
\int_{0}^1\int_{0}^1\f{s^{\all+1}t^{2n-3}}{(\delta +s+t^2)^{n+2} }\, dt \, ds\leq \f{C_n}{1-\all} \delta ^{\all-1}.
$$
\ele
\begin{proof}We estimate the integrals $I$ of the integrand by a covering of $[0,1]\times[0,1]$:

$(i)$ $\delta \leq t^2\leq s$.
$$
I\leq \int_{\delta }^1\int_{t=0}^{\sqrt s}\f{s^{\all+1}t^{2n-3}}{s^{n+2} }\, dt \, ds\leq\int_{\delta }^1s^{\all-2}\, ds\leq \f{1}{1-\all} \delta ^{\all-1}.
$$

$(ii)$ $\delta \leq s\leq t^2$.
$$
I\leq \int_{\sqrt \delta }^1\int_{s=0}^{t^2}\f{s^{\all+1}t^{2n-3}}{t^{2n+4} } \, ds\, dt\leq\int_{\sqrt \delta }^1t^{2\all-3}\, dt\leq \f{1}{1-\all} \delta ^{\all-1}.
$$

$(iii)$ $t^2\leq \delta \leq  s$.
$$
I\leq \int_{  0}^{\sqrt \delta }\int_{s=\delta }^{1} \f{s^{\all+1}t^{2n-3}}{s^{n+2} }\, dt \, ds= \delta ^{\all-n}\delta ^{n-1} =  \delta ^{\all-1}.
$$

$(iv)$ $ s\leq \delta \leq t^2$.
$$
I\leq \int_{0}^\delta \int_{t=\sqrt \delta }^{1}\f{s^{\all+1}t^{2n-3}}{t^{2n+4} }\, dt \, ds\leq \delta ^{\all+2}\delta ^{-3}=  \delta ^{\all-1}.
$$

$(v)$  $t^2\leq s\leq \delta $.
$$
I\leq \int_{0}^\delta \int_{t=0}^{\sqrt s}\f{s^{\all+1}t^{2n-3}}{\delta ^{n+2} }\, ds \, dt\leq \delta ^{-n-2}\int_{0 }^\del s^{n+\all}\, ds\leq  \delta ^{\all-1}.
$$

$(vi)$ $s \leq t^2\leq \del$.
$$
I\leq \int_{0}^{\sqrt\del}\int_{s=0}^{t^2}\f{s^{\all+1}t^{2n-3}}{\del^{n+2} }\, ds \, dt \leq   \delta ^{\all-1}.\qedhere
$$
\end{proof}

\pr{full-}Let  $r>1$. Let $D,\Phi, g^1$ be as in \rpa{estK}.  Suppose that $f\in C^1(\cc^n)$ is a function vanishing in $D$.  Then
\gan
\|H_0f \|_{\Lambda_r(\ov D)}\leq C_r(|W|_1)\| f\|_{\Lambda_r( {\cL U} \setminus D)}, \quad r>1,\\
|\pd_z^2H_0f(z)|\leq C_1(|W|_1)\dist(z,\pd D)^{-1}|f|_{1,\cL U\setminus D}, \quad z\in D.
\end{gather*}
\epr
\begin{proof}
Let $k=[r]\geq1$. We first consider the case $f\in C^r(\ov D)$. The above proof for $H_if$ with $i>0$ can be adapted easily. Let $\pd_z^{k+1}H_0f$ be a $(k+1)$-th order derivative of $H_0f$. It is a linear combination of
$$
Kf(z)=\int_{{\cL U} \setminus\ov  D}f(\zeta)\pd_{z}^{k+1}\left\{\f{A( W_1(z,\zeta), z,\zeta )}
 {\Phi^n( z,\zeta )}\right\}\, dV(\zeta).
$$
Let $z_0\in\pd D$.
Using a partition of unity, we may assume that for a neighborhood $B_0$ of $z_0$ in $\cc^n$ and for some $j$, we have
\ga
 \nonumber 
 \supp f\subset B_0\setminus D; \quad
 u( z,\zeta ):=\pd_{\zeta_j}\Phi( z,\zeta )\neq0, \quad  z,\zeta \in B_0.
\nonumber 
\end{gather}
Applying integration by parts $k-1$ times, we write $Kf$ as a linear combination of $K_1f$ with
 $$
  K_1  f(x):=\int_{\cc^n\setminus\ov  D} \f{f_1( z,\zeta )}{\Phi^{n+2}( z,\zeta )}\, dV(\zeta),
  \quad\forall z\in   D$$
 with
$
f_1( z,\zeta )=A \bigl(W_1(z,\zeta), z,\zeta \bigr)\,\pd_\zeta^{\nu_\ell}W_1\cdots \pd_\zeta^{\nu_1}W_1\,
\pd_\zeta^{\nu_{0}}f
$ and
\ga
\nonumber 
 \nu_0+\dots+\nu_\ell=k-1.
\end{gather}
Since $f(\zeta)=0$ in $\ov D$, it is easy to see that $  K_1  f\in \cL C^\infty(D)$.
We have for $\all>0$
$$
|K_1f(z)|\leq C(|W|_1)|f|_r \int_{s=0}^1\int_{t=0}^1\f{s^{\all+1}t^{2n-3}}{( d(z)+s+t^2)^{n+2}}\, ds\,dt\leq
C'(|W|_1)|f|_r  d(z)^{\all-1}.
$$
This gives us the desired estimate when $r$ is non integer.   When $r$ is a positive integer, the estimate follows from interpolation by \rp{interpolation}.
\end{proof}

\setcounter{thm}{0}\setcounter{equation}{0}
\section{Regularized Henkin-Ram{\'{\i}}rez functions}
\label{sect:5}

We now discuss our result for strictly pseudoconvex domains.
We first strengthen the classical  Henkin-Ram{\'{\i}}rez functions via the following result.
\pr{hrfn}
 Let $D$ be a bounded  strictly pseudoconvex domain in $\cc^n$ with $C^2$ boundary. Suppose that $\rho_0\in C^{2}({\cL U}  )$, $ \pd D=\{z\in {\cL U} \colon \rho_0=0\}$,  and $\pd \rho_0\neq0$ in $\pd D$.
Let $D_\del=\{z\in\cc^n\colon \dist(z,D)<\del\}$ and $D_{-\del}=\{z\in D\colon \dist(z,\pd D)
  >\del\}$. Let $\rho=E_2(e^{L_0\rho_0}-1)$ be a regularized $C^2$ defining function of $D$, where $L_0>0$ is sufficiently large so that $e^{L_0\rho_0}-1$ is  strictly plurisubharmonic in a neighborhood $\om$ of $\pd D$.
 There exist  $\del>0$ and functions $W$ satisfying the following.
 \bppp
 \item $W$ is defined in $D_{\del}\times (D_{\del}\setminus D_{-\del})$, $\Phi(z,\zeta)=W(z,\zeta)\cdot(\zeta-z)\neq0$ for $\rho(z)\leq \rho(\zeta)$ and $\zeta\neq z$, $W(\cdot,\zeta)$ is holomorphic in $D_\del$ for $z\in D_\del$, and $W\in C^{ 1}(D_\del\times ( D_{\del}\setminus D_{-\del}))$.
\item If    $|\zeta-z|<\epsilon$ and $\zeta\in D_\del\setminus D_{-\del}$, then
$\Phi(z,\zeta)=F(z,\zeta)
M(z,\zeta)$, $M(z,\zeta)\neq0$ and
  \ga{}
  \nonumber %
   F(z,\zeta)
  =-\sum\DD{ r }{\zeta_j}(z_j-\zeta_j)+ \sum a_{jk}(\zeta)(z_j-\zeta_j)(z_k-\zeta_k),\\
\nonumber 
  \RE F(z,\zeta)\geq  \rho (\zeta)- \rho (z)+|\zeta-z|^2/C, 
  \end{gather}
  with $(M,F) \in C^{1}( {D_{\del}}\times({D_{\delta}}\setminus D_{-\del}))$ and $a_{jk}\in C^{\infty}(\cc^n)$.
\item
For $ z\in D_\del$ and $
\zeta\in D_{\del}\setminus \ov D$, we have
\eq{chiW}
|\pd_z^i\pd_\zeta^jW(z,\zeta)|\leq C_{i,j}(D,L_0,|\rho_0|_{\ov D,2}) \sum_{j_1+j_2=j}\del^{-i-j_1}\Bigl\{1+d(\zeta,\pd D)^{1-j_2}\Bigr\}.
\eeq
 \eppp
The $(W_1,\dots, W_n)$ is called  a {\it regularized}  Henkin-Ram{\'{\i}}rez map.
  \epr
  \begin{proof} When $\rho$ is strictly plurisubharmonic,  the proof for (i) and (ii)
is in  {\O}vrelid~\ci{Ov71}  and see also Henkin-Leiterer
\cite{HL84}*{Thm.~2.4.3, p.~78;   Thm.~2.5.5, p.~81}, and Range~\cite{Ra86}*{Prop.~3.1, p.~284}.
In~\ci{GK14}, the Henkin-Ram{\'{\i}}rez functions for a family of strictly pseudoconvex domains are studied.
  Therefore,  only (iii) is new.


Fix $\del_0$ so that $\ov{D_{\del_0}}\setminus D$ is contained in $\cL U\cap\om$.  We have
$$\del_1=\min
\{\rho(\zeta)\colon \zeta \in\pd D_{\del_0}\}>\del_0/C.
$$
We have
$$
\sum_{j,k}
\DD{^2\rho(\zeta)}{\zeta_j\ov\pd \zeta_k}t_j\ov t_k\geq C_0|t|^2, \quad \zeta\in \om,
$$
with $C_0>0$. Define
$
D_{c}^*=\{z\in D_{\del_0}\colon \rho(z)<c\}.
$
We take
\eq{Fzze}
F(z,\zeta):=-\sum\DD{ \rho }{\zeta_j}(z_j-\zeta_j)- \sum a_{ij}(\zeta)(z_i-\zeta_i)( z_j-\zeta_j),
\eeq
where $a_{jk}\in C^{\infty}(\cc^n)$ with $|a_{jk}(\zeta)-\DD{^2}{\zeta_j\zeta_k}\rho|<1/C$ for $\zeta\in \cL U$.

 Fix $\e$ sufficiently small  so that for $|\zeta-z|<\e$ and $\zeta,z\in D_{\del_0}\setminus D_{-\del_0}$,
$$
\RE F(z,\zeta)\geq r(\zeta)-r(z)+|z-\zeta|^2/C_0.
$$
Let $\chi$ be a $C^\infty$ function satisfying $\chi(\zeta)=1$ for $|\zeta|<{3\e}/{4}$ and $\chi(\zeta)=0$ for
 $|\zeta|>{7\e}/{8}$.  Take $\del_2<\f{1}{4C_0}(\f{3}{4}\e)^2$ and $\del_2\in(0,\del_1)$.
For $z \in D_{\del_2}^*$, $\zeta\in D_{\del_1}^*\setminus D_{-\del_2}^*$,  and $|\zeta-z|>3\e/4$,
we have
$$
\RE F(z,\zeta)\geq r(\zeta)-r(z)+|z-\zeta|^2/C_0\geq -2\del_2+\f{1}{C_0}\left(\f{3}{4}\e\right)^2>\f{1}{2C_0}\left(\f{3}{4}\e\right)^2.
$$
Thus we can define
\eq{fzztwo}
f(z,\zeta)=\begin{cases} \db_z(\chi (\zeta-z)\log F(z,\zeta))
      & \text{if $3\e /4<|\zeta-z|<7\e/8$}, \\
   0   & \text{otherwise},
\end{cases}
\eeq
for  $z \in D_{\del_2}^*$ and $\zeta\in D_{\del_1}^*\setminus D_{-\del_2}^*$.
   Define
\ga\label{defuz}
u(z,\zeta)=T_{D^*_{\del_2}}f(\cdot,\zeta)(z), \quad\forall z \in D_{\del_2}^*,\
\zeta\in D_{\del_1}^*\setminus D_{-\del_2}^*.
\end{gather} Here   $T_0=T_{D_{\del_2}^*}$ is a linear $\db$ solution operator
that  admits an interior super-norm estimate on $D_{\del_2}^*$ (see~\cite{HL84}*{Thm.~2.3.5, p.~76}). Namely, for any $\db$-closed $(0,1)$-form $\var$ in $D_{\del_2}^*$, $\db T_0\var=\var$ and
\eq{ukT0}
|T_0\var|_{C^0(D^*_{\del_3})}\leq C_0^*|\var|_{C^0(D_{\del_2}^*)},
\eeq
for $\del_3\in(0,\del_2)$.
By the linearity of $T_0$ and \rp{estL},  the $ u(z,\zeta)$ and $\pd_\zeta u(z,\zeta)$ are uniformly continuous in $D^*_{\del_2}\times( D^*_{\del_1}\setminus D^*_{-\del_2})$.  We also have
$$ 
|f(z,\cdot)|_{ D^*_{\del_1}\setminus D_{-\del_2}^*;1}\leq C(|\rho|_2),\  z\in D_{\del_2}^*;
\quad |u(z,\cdot)|_{ D^*_{\del_1}\setminus D_{-\del_2}^*,1}\leq C_0^*C(|\rho|_{2}), \  z\in D^*_{\del_3}.
$$ 
Define for $z\in D^*_{\del_3}$ and $\zeta\in D_{\del_1}^*\setminus D_{-\del_2}^*$,
\eq{defPhi}
\Phi(z,\zeta):=\begin{cases}
 F(z,\zeta)e^{-u(z,\zeta)}    & \text{if $|\zeta-z|\leq 3\e/4$}, \\
e^{\chi(\zeta-z)\log F(z,\zeta)-u(z,\zeta)}    & \text{otherwise}.
\end{cases}
\eeq
Then
$
|\Phi(\cdot,\zeta)|_{D^*_{\del_3},0}\leq C|\rho|_{2}$ for $
\zeta\in D_{\del_1}^*\setminus D_{-\del_2}^*$. Also
$$
|\Phi(z,\cdot)|_{D_{\del_1}^*\setminus D_{-\del_2}^*,1}\leq C(C_*,|r|_{2}), \quad z\in D_{\del_3}^*.
$$
Fix $\del_4^*\in(0,\del_3^*)$. By Hefer's decomposition theorem~\cite{HL84}*{p.~81}, there are continuous linear mappings
$
T_j\colon \cL O(D^*_{\del_3})\to \cL O((D^*_{\del_4})^2)
$
so that
$
h(\tilde z)-h(z)=\sum_{j=1}^nT_jh(z,\tilde z)(\tilde z_j-z_j).
$
   Then we have
\eq{defnTjv}
\Phi(\tilde z,\zeta)-\Phi( z,\zeta )=\sum T_j\Phi(\cdot,\zeta)(z,\tilde z)(\tilde z_j-z_j).
\eeq
Set $W_j( z,\zeta )=T_j\Phi(\cdot,\zeta)( z,\zeta ).
$
We know that $T_j\Phi(\cdot,\zeta)( z,\eta)$ is holomorphic in $z,\eta$.  We express the boundedness of $T_j$ as
\eq{ukTj}
 |T_jh|_{C^0((D^*_{\del_4^*})^2)}\leq C_j^*|h|_{C^0(D^*_{\del_3})},\quad h\in \cL O(D^*_{\del_3}),
 \quad j=1,\dots, n.
\eeq
The linearity and continuity of $T_j$ imply that $W_j$ and its first-order derivatives in $\zeta$
are continuous. Since $W_j$ is holomorphic in $z$,   the Cauchy formula implies that $W$ is in $C^1(D_\del\times (D_\del\setminus D_{-\del}))$ by shrinking $\del$ slightly.

We now use the fact that $\rho$ is a regularized $C^2$ defining function for the domain $D$ to   estimate  the higher order derivatives of $\Phi(z,\zeta)$ for $\zeta\not\in\ov D$. We   restrict
$$
z\in D_{\del_5^*}, \quad \zeta\in D^*_{\del_4}
\setminus \ov D.
$$
Here we take $\del_5^*\in(0,\del_4^*)$. This also allows us to use Cauchy inequality in the $z$ variables.

By \re{fzztwo}, \re{defuz} and the linear estimate \re{ukT0}, we first   see that for each $j$,
 $\pd_\zeta^j u(z,\zeta)$
are continuous in $(z,\zeta)\in D_{\del_3}^*\times(D_{\del_1}^*\setminus \ov D)$. Moreover,
$$
|\pd_\zeta^ju(\cdot,\zeta)|_{D_{\del_3}^*;0}\leq C_j(C_0^*,|\rho|_2)(1+d(\zeta)^{1-j}), \quad \zeta\in D_{\del_1}^*\setminus \ov D.
$$
Here we have use $|\pd_\zeta^j\rho(\zeta)|\leq C_j(1+d(\zeta)^{1-j})$ as well as the product rule for
$$
\log F(z,\zeta)=\log \RE F(z,\zeta)+\log\Bigl(1+i\f{\IM F(z,\zeta)}{\RE F(z,\zeta)}\Bigr),\quad 3\e/4<|\zeta-z|<7\e/8,
$$
where $z\in D^*_{\del_3},\zeta\in D_{\del_1}^*\setminus \ov D$.
By \re{defPhi},  we get $|\pd_\zeta^j\Phi(\cdot,\zeta)|_{D_{\del_3}^*;0}\leq C_j(1+d(\zeta)^{1-j})$. Here and in what follows, we let
$$
C_j:= C_j(C_0^*,\dots, C_n^*,|\rho|_2).
$$
By the linearity and continuity of $T_j$ and the holomorphicity of $T_j\Phi(\cdot,\zeta)( z,\eta)$ in $\eta$, we have
$$ 
\pd_\zeta^\all W_\ell( z,\zeta )=\sum\binom{\all}{\beta}\pd_\eta^{\all-\beta}\Bigr|_{\eta=\zeta}\bigr.
T_\ell\pd_\zeta^\beta\Phi(\cdot,\zeta)(z,\eta)
$$ 
for $z\in D_{\del_4}^*$ and $\zeta\in D_{\del_4}^*\setminus \ov D$. By the linearity of estimate \re{ukTj} for $T_\ell$ and Cauchy inequalities applied to the last term, we get
 $$|\pd_\zeta^jW(\cdot,\zeta)|_{D_{\del_5}^*;0}\leq C_j \sum_{j_1+j_2=j}\dist(D^*_{\del_5},\pd D_{\del_4}^*)^{-j_1}
 (1+d(\zeta)^{1-j_2})$$ for $j=1,2,\dots.$
   By Cauchy inequalities, we get
 $$
|\pd_\zeta^jW(\cdot,\zeta)|_{D_{\del_5}^*;i}\leq C_j \sum_{j_1+j_2=j}\dist(D^*_{\del_5},\pd D_{\del_4}^*)^{-i-j_1}
 (1+d(\zeta)^{1-j_2})
 $$
for $\zeta\in D_{\del_5}^*\setminus \ov D$.  Finally,  we fix $\del\in(0,\del_5)$. We have achieved
  \re{chiW}.
\end{proof}

  \th{full}
Let $D=\{z\in  {\cL U} \colon \rho_0<0\}$ be a strictly pseudoconvex domain with $C^2$ boundary that is relatively compact in $\cL U$, where $\rho_0\in C^2(\cL U)$ and $d\rho_0\neq0$ in $\pd D$.  Let $H_q$ be defined by \rea{Hqv} and \rea{H0f}, where $g^1=W$ is the regularized Henkin-Ram{\'{\i}}rez function $D_\del\times (D_\del\setminus D_{-\del})$ as in \rpa{hrfn}
 and $\Phi( z,\zeta )=W(z,\zeta)\cdot (\zeta-z)$. Let $\var$  be a $(0,q)$-form   such that $\var, \db\var$ are in $C^1(\ov D)$. Then in $D$
\ga\label{tsqfG}
\var=\db H_q\var+H_{q+1}\db\var, \quad 1\leq q\leq n,\\
\var_{0}=H_0\var_{0}+H_1\db \var_{0}.
\label{HqvG}
\end{gather}
Moreover,  we have
\ga
\label{hqv.5}\|H_q\var  \|_{\Lambda_{r+1/2}(\ov D)}\leq C_r(D) \|\var\|_{\Lambda_r(\ov D)},\quad r>1,  q>0,\\
\label{hqv.51}\|H_q\var  \|_{D;3/2}\leq C_1(D)\|\var\|_{D;1}, \quad q>0,\\
\|H_0\var  \|_{\Lambda_r(\ov D)}\leq C_r(D)\| \var \|_{\Lambda_r(\ov D)}, \quad r>1,\label{hqv.52}
\\
\label{hqv.53}\|\pd_z^2H_0\var(z)  \| \leq C_1(D)\dist(z,\pd D)^{-1}\|\var \|_{D;1}, \quad z\in D.
\end{gather}
The constants $C_1(D), C_r(D)$   are stable under  a small $C^2$ perturbation. They depend on the $\e, N, M$ in \rda{mins}, the $L$ in \rla{estR}, $|\pd\rho_0|_{D\setminus D_{-\del_0};0}$, $|\pd\db\rho_0|_{D\setminus D_{-\del_0};0}$,
$|\pd \rho_0|_{\ov D;2}$, as well as the constants $L_0, \del, C_0^*,\dots, C_n^*$ in the proof of \rpa{hrfn}.
Therefore, $C_r(\widetilde D)\leq  \widehat C_r(D,\e)<\infty$ for all $r\in(0,\infty)$
when $\widetilde D$ has a defining function $\tilde \rho$ such that $|\tilde \rho_0-\rho_0|_{\cL U;2}<\e$ for sufficiently small $\e$.
\end{thm}
\begin{proof}Let $d(z)=\dist(z,\pd D)$. Let $\rho=E_2(e^{L_0\rho_0}-1)$ as in \rp{hrfn}.   Let us first choose the local coordinates described in \rp{estK}. As in \ci{HL84}*{p.~73}, we take
$$
s_1=\rho(\zeta), \quad s_2=\IM(\rho_\zeta\cdot(\zeta-z)), \quad t=(\RE(\zeta'-z'),\IM(\zeta'-z')).
$$
Let $F$ be as in \rp{hrfn}. First, we have $|F(z,\zeta)|\geq\RE(F(z,\zeta))\geq c_*|\zeta-z|^2$.  Note that $\rho(\zeta)\approx d(\zeta)$ and $-\rho(z)\approx d(z)$.  Then we have $2|F(\zeta,z)|\geq |s_2|+\rho(\zeta)-\rho(z)+c_*|\zeta-z|^2\geq c_*'(d(z)+s_1+|s_2|+ |t|^2).$ Since $z\in D$ and $\zeta\not\in\ov D$, we also have $d(z)\leq|\zeta-z|, d(\zeta)\leq|\zeta-z|$, $\rho(\zeta)\approx d(\zeta)\leq|\zeta-z|$. Hence
$
d(z)+d(\zeta)+s_1+|s_2|+|t|\leq C|\zeta-z|.
$
Thus, we have   verified \re{LbPhi}-\re{LbPhi+} in which $\Phi=FM$. We obtain the desired   estimates by \rp{hrfn}, \rp{estL}, \rp{estK} (ii), and \rp{full-}.
\end{proof}

\setcounter{thm}{0}\setcounter{equation}{0}
\section{Boundary regularities  of the elliptic  differential complex for the Levi-flat Euclidean space}
\label{sect:6}

We consider the complex for the exterior differential
$\cL D:=d_{t}+\db_z$,  where  $(z,t)$ are coordinates of $\cc^n\times\rr^M$.  We will also write it as
$\cL D=d^0+\db.
$

The Poincar\'e lemma for a $q$-form on a bounded star-shaped domain $S$ has the form
\gan 
\nonumber
\phi=d^0R_q\phi+R_{q+1}d^0\phi, \quad q>0;\quad \phi=R_1d^0\phi+\phi(0), \quad q=0; \\
   R_q\phi(t):=\int_{\theta\in[0,1]}H^*\phi(t,\theta). \nonumber
\end{gather*}
Here $H(t,\theta)=\theta t$ for $(t,\theta)\in S\times[0,1]$; see~\cite{Sp99}*{p.~224} and
\cite{Tr92}*{p.~105}.
If  $\phi(t)=f(t)dt_1\wedge \cdots\wedge dt_q$, we have
$$
  R_q\phi(t)=\left\{\int_0^1f(\theta t)\theta^{q-1}d\theta\right\} \sum(-1)^{j-1} t_j\,dt_1\cdots\widehat{dt_j}\cdots  dt_{q}.
$$
It is immediate
that for $q>0$
\eq{3p13}
 |R_q\phi |_{S;r}\leq C_r  |\phi |_{S;r}, \quad 0\leq r<\infty.
\eeq
Here $C_r$ depends only on the diameter of $S$.
By the  interpolation argument, we obtain
\eq{3p13z+}
 |R_q\phi |_{\Lambda_r(\ov S)}\leq C_r(D)\rho |\phi |_{\Lambda_r(\ov S)}, \quad 0<r<\infty,
\eeq
provided $S$ is a bounded Lipschitz domain.

A differential form $\var$ is called of {\it mixed type} $(0,q)$
if
$
\var=\sum_{i=0}^ {q} [\var]_{i},
$
where
$$[\var]_{ i}=\sum_{|I|=i, |I|+|J|=q} a_{IJ}d\ov z^I\wedge dt^J.
$$  Thus
$
  [\var]_i=0$, for $i>n.
$
The $\cL D $ acts on a function $f$ and a $(0,q)$-form as follows
$$
\cL Df=\sum\f{\pd f}{\pd t_m}dt_m+\sum\f{\pd f}{\pd\ov z_\all}\, d\ov z_\all, \quad
\cL D \sum a_{IJ}d\ov z^I\wedge dt^J=\sum \cL Da_{IJ}\wedge d\ov z^I\wedge dt^J.$$
We have
$
\cL D ^2=0$ and $   d^0\db+\db d^0=0$.
We also have
\ga \nonumber 
[\cL D \var]_0=d^0[\var]_{0},\quad
[\cL D \var]_i=d^0[\var]_{i+1}+\db[\var]_{i}, \quad 0<i\leq n.
\end{gather}
For   $\var=\sum \var_{IJ}d\ov z^I\wedge dt^J=\sum\tilde\var_{IJ}dt^J\wedge d\ov z^I$ on $\ov D\times \ov S$, define
\gan
H_i\var=\sum_{|I|=i}H_i(\var_{IJ}d\ov z^I)\wedge dt^J, \quad
R_{i}\var =\sum_{|J|=i}R_{i}( \tilde \var_{IJ}d t^J)\wedge d\ov z^I.
\end{gather*}
Thus $H_i\var=H_i[\var]_i$,  while $R_{q-i}\var=R_{q-i}[\var]_i$ if $\var$ has the (mixed) type $(0,q)$.
\begin{defn} Let
 $0< r\leq1$  and $0\leq\all<1$.  Let $\Lambda_*^{r,0}(\ov D\times\ov S)$
   be the set of continuous functions $f$ in $\ov D\times\ov S$ so that
$t\to |f(\cdot,t)|_{\Lambda_r(\ov D)}$  is bounded  in $\ov S$.  For $a>1$ and $k\in\nn$, let $\Lambda_*^{a,k}(\ov D\times\ov S)$ be the set of functions $f$ so that $\pd_z^i\pd_t^jf$ are in $ \Lambda_*^{a-i,0}(\ov D\times\ov S)$    for $j\leq k$ and $i<a$. Define $C_*^{a,k}(\ov D\times\ov S)$ analogously.
\end{defn}

We now derive the following homotopy formulae.
\pr{dpht} Let $1\leq q\leq n+m$. Let $D$ be a bounded strictly pseudoconvex domain in $\cc^n$ with $\pd D\in C^{2}$ and let $S$ be a bounded domain in $\rr^m$ so that $\theta S\subset S$ for $\theta\in[0,1]$.
  Let $\var$ be a mixed $(0,q)$-form in $ D\times S$.
\bppp
\item If $\var\in C^{1,1}_*(\ov D\times\ov S)$ and  $\db\var$ are in $C^{1,0}_*(\ov D\times\ov S)$, then
\ga\label{vdhv}
\var=\cL DT_q\var+T_{q+1}\cL D \var,\\
  T_q\var=R_qH_0[\var]_0+\sum_{i>0}H_{i}[\var]_i.
\label{hvar+}
\end{gather}
\item If  $\var\in C^{1}(\ov D\times\ov S)$ and $\db[\var]_q(\cdot,0)\in C^1(\ov D)$, then
\ga
\label{vdhv+}
\var=\cL D \widetilde T_q\var+\widetilde T_{q+1}\cL D \var,\\
\label{hvar}
\widetilde T_q\var(z,t)=H_q[\var]_q(\cdot,0)(z)+\sum_{i< q}R_{q-i}[\var]_i(z,\cdot)(t).
\end{gather}
\eppp
\epr
\begin{proof} For the related homotopy formulae when $H_i$'s are replaced by the Leray-Koppelman homotopy operators, see  Treves~\ci{Tr92}*{Sect.~VI.7.12,  Sect.~VI.7.13, p.~294} for suitable $q's$ and \ci{Go16} for   arbitrary $q$.


$(i)$ Recall that the homotopy operators $H_i$ are linear. To derive the homotopy formulae for $\cL D$, we will use the following estimates from \re{hqv.51}:
$$
|H_i[\var]_i|_{\ov D;0}\leq C|[\var]_i|_{\ov D;1}, \quad i=0,1,\dots n.
$$
Thus if $\psi_j$ converges to $\psi$  in $C^1(\ov D)$ norm, then
\eq{intlim}
\lim_{j\to\infty} H_i\psi_j=H_i\psi.
\eeq
Also for $\psi\in C^{1,1}_*(\ov D\times\ov S)$, we have
\eq{pdHpd}
\f{\pd}{\pd t_j}H_i\psi(\cdot,t)=H_i\f{\pd }{\pd t_j}\psi(\cdot,t).
\eeq
Analogously, if $\psi_j$ converges to $\psi$ in $C^0(\ov S)$,  then
$
\lim_{j\to\infty} R_i\psi_j=R_i\psi.
$
For $\psi\in C^{1,0}_*(\ov D\times\ov S)$, we have
$
\f{\pd}{\pd \ov z_j}R_i\psi=R_i\f{\pd}{\pd \ov z_j}\psi.
$
Note that $\db$ commutes with the pull-back $H^*$ of $H(t,\theta)=\theta t$.  By \re{checksign}, we have
\ga  \label{d0p}
d^0H_i\var=-H_id^0\var, \quad \var\in C^{1,1}_*(\ov D\times\ov S),\\
\db R_{i}\var=-R_{i}\db\var, \quad \var\in C_*^{1,0}(\ov D\times\ov S).
\label{d0p+}
\end{gather}

Let us start with the integral representation of $[\var]_0$ in $D$.  Since $\var$ has total degree $q$,
then $\deg_x[\var]_0=q>0$.
We apply \re{HqvG} for functions and the Poincar\'e formula for $d^0$ by \re{d0p}. Thus for $[\var]_0\in C^{1,1}_*$, we   obtain in $D\times S$
\eq{bmkp}
[\var]_0=H_0[\var]_0+H_1\db[\var]_0=(d^0R_qH_0[\var]_0+R_{q+1}d^0H_0[\var]_0)+H_1\db[\var]_0.
\eeq
Since $\db_z\Om_{0,0}^1( z,\zeta )=0$, we have
$
d^0R_qH_0[\var]_0=\cL DR_qH_0[\var]_0$.
Combining it with $d^0[\var]_0=[\cL D\var]_0$, we express \re{bmkp} as
 \eq{bmkp+}
[\var]_0= \cL DR_qH_0[\var]_0+R_{q+1}H_0[\cL D\var]_0+H_1\db [\var]_0.
\eeq

Analogously, for $[\var]_j\in C^{1,1}_*(\ov D\times\ov S)$, we get
\eq{varj1}
[\var]_j=\db H_j[\var]_j+H_{j+1}\db[\var]_j=\cL DH_j[\var]_j-d^0H_j[\var]_j+H_{j+1}\db[\var]_j.
\eeq
 By \re{d0p} and $d^0[\var]_j=[\cL D \var]_j-\db[\var]_{j-1}$, we obtain
\aln
\sum_{j>0}\left(-d^0H_j[\var]_j+H_{j+1}\db[\var]_j\right)&=\sum_{j>0}\left(H_j[\cL D \var]_j -H_{j}\db[\var]_{j-1}+H_{j+1}\db[\var]_j \right)\\
&= -H_{1}\db[\var]_{0}+\sum_{j>0}H_j[\cL D \var]_j.
\end{align*}
 Here we have used $H_{n+1}=0$.
Combining it with \re{bmkp+} and \re{varj1}, we obtain
\aln
\var=\cL DR_qH_0[\var]_0+R_{q+1}H_0[\cL D \var]_0+\sum_{j>0}\cL D  H_j[\var]_j+\sum_{j>0} H_j[\cL D \var]_j,
\end{align*}
which gives us (i).


$(ii)$ By $[\var]_j\in C^1(\ov S\times\ov S)$ and
  the Poincar\'e lemma, we obtain
\aln
\var&=[\var]_q+\sum_{i<q}(d^0R_{q-i}[\var]_i+R_{q+1-i}d^0[\var]_i)=[\var]_q+\sum_{i<q}\cL DR_{q-i}[\var]_i+R_{q+1-i}\cL D[\var]_i.\end{align*}
Here we have  used $\db R_{q-i}[\var]_i=-R_{q-i}\db[\var]_i$ for $i<q$ by \re{d0p} and
$
R_{q+1-i}d^0[\var]_i=R_{q+1-i}\cL D[\var]_i-R_{q-i}\db[\var]_i.
$
 We express
$$
\sum_{i<q}R_{q+1-i}\cL D[\var]_i=\sum_{i<q}R_{q+1-i}([d^0\var]_i+[\db\var]_{i+1})= -R_1[d^0\var]_{q}+\sum_{i\leq q}R_{q+1-i}[\cL D\var]_{i+1},
$$
because $[\db\var]_0=0$. We have
 $
 [\var]_q(z,t)-R_{1}d^0[\var]_q(z,\cdot)(t)=[\var]_q(z,0).
 $
 We now apply the  homotopy formula \re{tsqf} to obtain
 $$
 [\var]_q(\cdot,0)=\db H_q[\var]_q(\cdot,0)+H_{q+1}\db[\var]_q(\cdot,0)=\cL DH_{q}[\var]_q(\cdot,0)+H_{q+1}([\cL D\var]_{q+1}(\cdot, 0)).
  $$
Combining the identities, we get $(ii)$. 
\end{proof}

\th{bD-reg} Let $q>0$.
  Let $D$ be a strictly pseudoconvex domain with $C^{2}$ boundary. Let  $S$ be a bounded star-shaped domain in $\rr^m$. Let $\var$ be a $\cL D$-closed $(0,q)$-form in $C^1(\ov D\times \ov S)$. Then there exists a solution $u\in C^{1}(\ov D\times \ov S)$ to $\cL Du=\var$. Furthermore, the following properties hold.
\bppp
     \item Suppose that $[\var]_0=0$ and $\db\var\in C^{1,0}_*(\ov D\times\ov S)$.
         If   $\var\in\Lambda^{r,k}_*(\ov D\times\ov S)$ with  $k\in\{0,1,\dots,\infty\}$  and $r\in(1,\infty]$, the $u$ is in $\Lambda_*^{r+1/2,k}(\ov D\times\ov S)$.
\item Let $r\in[1,\infty]$. If
  $\var\in C^{r}(\ov D\times\ov S)$,  the $u$ is in $ C^{r}(\ov D\times\ov S)$.
  If $\var\in \Lambda_{r}(\ov D\times\ov S)$ and $S$ is a Lipschitz domain,  the $u$ is in $ \Lambda_{r}(\ov D\times\ov S)$ for $r>1$.
  \eppp
\eth
\begin{proof}
$(i)$ follows from \re{hvar+} and   \re{hqv.5}.
$(ii)$ follows from \re{hvar}, \re{3p13}, and \re{hqv.52}. Indeed, we have $\db[\var]_q(\cdot,0)=0$ as  $\cL D\var=0$.
We first obtain the assertion when $r$ is non-integer. The general case is obtained via interpolation for the  Lipschitz domain $D\times S$.
\end{proof}


\newcommand{\doi}[1]{\href{http://dx.doi.org/#1}{#1}}
\newcommand{\arxiv}[1]{\href{https://arxiv.org/pdf/#1}{arXiv:#1}}

  \def\MR#1{\relax\ifhmode\unskip\spacefactor3000 \space\fi%
  \href{http://www.ams.org/mathscinet-getitem?mr=#1}{MR#1}}


\end{document}